\documentclass[preprint,12pt]{elsarticle}

\usepackage{amssymb}
\usepackage{amsthm}
\usepackage{amsmath}
\usepackage{color}
\usepackage[labelfont=bf,list=true,font=footnotesize]{subcaption}

\journal{Elsevier}

\begin{document}

\begin{frontmatter}



\title{A nonconservative macroscopic traffic flow model in a two-dimensional urban-porous city.}


\author[guadalajara]{N. Garc\'\i a-Chan\corref{cor1}{}}
\ead{nestorg.chan@cucei.udg.mx}
\address[guadalajara]{Depto. F\'\i sica. C.U.~Ciencias Exactas e Ingenier\'\i as.
Universidad de Guadalajara.\\ 44420 Guadalajara. Mexico}
\author[vigo]{L.J. Alvarez-V\'azquez}
\cortext[cor1]{Corresponding author}
\ead{lino@dma.uvigo.es}
\address[vigo]{Depto.~Matem\'atica Aplicada II. E.I.~Telecomunicaci\'on.
Universidad de Vigo. \\ 36310 Vigo. Spain}
\author[vigo]{A. Mart\'\i nez}
\ead{aurea@dma.uvigo.es}
\author[lugo]{M.E.~V\'azquez-M\'endez}
\ead{miguelernesto.vazquez@usc.es}
\address[lugo]{Depto. Matem\'atica Aplicada. E.P.S. Universidad de
Santiago de Compostela. \\ 27002 Lugo. Spain}

\begin{abstract}
In this paper we propose a novel traffic flow model based on understanding the city as a porous media, this is, streets and building-blocks characterizing the urban landscape are seen now as the fluid-phase and the solid-phase of a porous media, respectively. Moreover, based in the interchange of mass in the porous media models, we can model the interchange of cars between streets and off-street parking-spaces. Therefore, our model is not a standard conservation law, being formulated as the coupling of a non-stationary convection-diffusion-reaction PDE with a Darcy-Brinkman-Forch\-heimer PDE system. To solve this model, the classical Galerkin $P_1$ finite element method combined with an explicit time marching scheme of strong stability-preserving type was enough to stabilize our numerical solutions. Numerical experiences on an urban-porous domain inspired by the city of Guadalajara (Mexico) allow us to simulate the influence of the porosity terms on the traffic speed, the traffic flow at rush-valley hours, and the streets congestions due to the lack of parking spaces.     

\end{abstract}




\end{frontmatter}


\section{Introduction}
One of the most worrying problems related to the urban development is the air pollution due to traffic flow. Not only the enormous number of cars (reaching hundreds of thousands, even millions in bigger cities), but also traffic jams, congested streets and poor planning turn metropolises in highly polluted places. To draw public strategies which mitigate the air pollution, {\it a priori} estimates of traffic flow on the entire city are needed, but the complexity of the traffic dynamics (as the cars move to their destinations on the countless streets of the city road network) along with their constant influx from residential zones and departure to parking-spaces make these a hard task. To our knowledge, the mathematical models based on partial differential equations (PDE) and their numerical solution could be an effective and cheap first-approximation in order to overcome this issue with, at least, a rough or mean estimation.  However, choosing an adequate and suitable mathematical model is not trivial, and new ideas are necessary. 

An important issue in the context of traffic flow modelling is doubtless the capture of the urban landscape with its streets, building-blocks, malls, schools, parks and other urban elements. In this sense, an interpretation of the urban landscape as a porous media allows us to use mathematical models for compressible and incompressible flows on porous media (see, for instance,  \cite{das2018_book}).  These models are formed by the equations of continuity for the fluid mass, the Darcy-Brinkman-Forchheimer equations for the fluid velocity, and advection-diffusion-reaction equations for the interchange of mass and energy.  In fact, the idea of an urban porous media to development numerical studies is not new.  For instance, different models for an incompressible flow on a porous city were recently applied in \cite{hu2012,wang2021,ming2021} to study numerically the intensity of the urban heat island (UHI) due to the different levels of  anthropogenic heat flux. Also, in \cite{garcia-chan2023} a numerical study of the UHI with a model formulated with the equations of mass, momentum and energy interchange on a porous-urban domain was capable to reproduce two well-known phenomenons: the inverse temperatures along the 24 hours of the day, and the transport of air heated by the wind field. 

Within the many traffic flow models existing in the literature (see \cite{treiber2013} and the references therein) we focus on the so-called macroscopic dynamic traffic assignment (DTA) models. These models are characterized by using a system of PDE integrated by the continuity and momentum equations, and its substantial motivation is based on the fact that traffic flow is assumed to have a similar behavior to a compressible gas flow. More frequently, these models are defined on one-dimensional domains as avenues or streets, but exceptionally two-dimensional domains are also considered, which opens the possibility of covering an entire city. This is possible only with the assumption that a high density of streets allows drivers to follow a straight direction towards its destination as a mean direction of the zigzag path on the real cities. In \citep{Jiang2014,Jiang2016} the authors use a macroscopic DTA model to simulate the traffic flow on a two-dimensional city, analysing its consequences on different pollutants emissions with a pseudo-conservation law of cars, treated as an inviscid flow. However, these papers do not consider the urban porosity in the sense of the proportion of the urban area covered by building-blocks and the complementary covered by streets, either consider porosity related to parking-spaces, with all cars leaving the city through a suitable boundary condition.
       
To our knowledge,  porous media models have not been previously applied to simulate the traffic flow in an entire city. Thus, in this paper we propose a novel model to deal with this unaddressed topic.  The outline of the paper is as follows. In section 2 we make a detailed deduction of our model, starting from recalling the Euler equations to model a compressible flow of inviscid fluid, passing for the macroscopic DTA traffic models which add to the Euler system a source of mass, diffusion effects, a relaxation term, and drop the pressure by a state relation, until reaching the macroscopic DTA corrected pressure model. This study of preceding models allows us to observe elements not considered previously, and formulate our model in a consistent manner with the idea of the porous city. In section 3 the process to derive the numerical scheme is displayed,  based on combining a finite element method of $P_1$ type with an explicit time marching scheme of strong stability preserving (SSP) type. In section 4 several numerical experiences with the city of Guadalajara as domain are presented, with the aim of evaluating the influence on traffic density and speed of the different key parameters of the model: porosity, relaxation time, rate of cars parking, and traffic sources. Finally, in section 5 some discussions and conclusions are written, highlighting the novelties of our model and results with respect to other authors.

\section{The mathematical model}
\label{sec:mathmodel}

\subsection{Traffic flow models as a conservation law and its pressure correction}
The well-known Euler equations were defined to model the compressible flow of a fluid without diffusion effects (i.e., inviscid). These equations can be formulated in primitive variables as follows \cite{huber2015}:
\begin{subequations}\label{ModEuler}
\begin{align}
\dot{\rho} + \nabla\cdot(\rho\mathbf{u}) = 0,\label{ModEuler:massconservation}\\
\dot{\mathbf{u}} + (\mathbf{u}\cdot\nabla)\mathbf{u} +\frac{\nabla P}{\rho} =\mathbf{0},\label{ModEuler:momentum} \\
\dot{P} + \mathbf{u}\cdot\nabla P + \rho c^2\nabla\cdot\mathbf{u}=0,\label{ModEuler:pressure}
\end{align}
\end{subequations}
where $\rho$ is the density, $\mathbf{u}=(u_1,u_2)$ is the speed, $P$ is the pressure, $c^2$ is the sonic speed, and the dot ($\,\dot{\,}\,$) will denote from now on the time derivative. The system is formed by the equation for mass conservation (\ref{ModEuler:massconservation}),  the equation for momentum conservation (\ref{ModEuler:momentum}), and the pressure equation in accordance with density and sonic speed (\ref{ModEuler:pressure}).  
The Euler system is hyperbolic, and can produce shock-waves and isotropic flows. This is of special interest to model phenomena in fluid dynamics \cite{huber2015} where the fluid density $\rho$ could present variations in space and time along with subsonic and ultrasonic speeds (compressible flow). Typical benchmarks are the oblique shocks, the reflected shocks, and the multi-phase flows \cite{tezduyar2006,caltagirone2011}. 

Nevertheless, the traffic flow is far to reach ultrasonic velocities,  do not present a multi-phase nature, and above mentioned shock waves benchmarks are not suitable for the traffic behaviour. Also, the pressure due to traffic has a straightforward relation to density, and the conservation of mass and momentum could be desired in particular cases but is not indispensable \cite{treiber2013,Jiang2016}. Therefore, the original Euler system is not appropriate to model the traffic flow, thus several alternatives have been developed as the so-called macroscopic with dynamical velocity models or the macroscopic dynamic assignment models. Those models assume that traffic flow behavior is similar to a compressible flow of gas, but with a pseudo-conservation of mass and momentum, and neglecting the pressure equation by a state relation in order to warranty an increasing relation with traffic density (see \cite{treiber2013} and references therein). Models like Payne-Whitham \cite{payne1971}, Keyner-Konh\"uauser \cite{kerner_konhauser1993}, and gas-kinetic \cite{herman-tenny-prigrione1971} include diffusion effects, viscosity, traffic demand (mass source), and a relaxation term with respect to a desired velocity. Those models could be formulated in the following generalized form:
\begin{subequations}\label{ModelTraffic}
\begin{align}
\dot{\rho} + \nabla\cdot (\rho\mathbf{u}) - \nabla\cdot(\nu\nabla\rho) = q,\label{MassConservation}\\
\dot{\mathbf{u}} +( \mathbf{u}\cdot\nabla)\mathbf{u} = -\frac{\nabla P}{\rho} + \frac{\mathbf{v}-\mathbf{u}}{\tau} + (\nabla\cdot\mu\nabla)\mathbf{u},\\
P = P(\rho),
\end{align}  
\end{subequations}
where $\nu$ is the diffusion coefficient of the traffic density, $\mu$ is the fluid viscosity, $\tau$ is the relaxation time, $\mathbf{v}=\mathbf{v}(\rho)$ is the desired speed, $q$ is the traffic demand, and $P=P(\rho)$ is a state relation of the pressure to the traffic density.

However, this kind of models (\ref{ModelTraffic}) present an important drawback which was pointed out in \cite{daganzo1995}: in one-dimensional macroscopic DTA models drivers will respond to stimuli from the front and from behind but, if fact,  drivers respond only to front stimuli. This inconsistency with the reality was corrected in \cite{Aw2000} with a suitable relation $P=P(\rho)$, and with the convective derivative $(\, \dot{ } + \mathbf{u}\cdot\nabla)$ applied to the sum of speed and pressure $(\mathbf{u}+P\mathbf{1})$ instead of only to $\mathbf{u}$, where $\mathbf{1}$ denotes the all-ones vector. The corrected pressure model proposed by \cite{Aw2000}, for the particular admissible choice $P(\rho)=c^2 \rho$,  is written here in its two-dimensional version as the following system:
\begin{subequations}\label{ModAx}
\begin{align}
\dot{\rho} + \nabla\cdot(\rho\mathbf{u}) = 0,\label{ModAx:massconservation}\\
\dot{\mathbf{u}} + (\mathbf{u}\cdot\nabla)\mathbf{u} - c^2\rho(\nabla\cdot\mathbf{u}) \mathbf{1}=\mathbf{0}. \label{ModAx:momentum}
\end{align}
\end{subequations} 

In \cite{Aw2000} the authors conclude that this new model presents good mathemat\-ical properties: is hyperbolic, has a unique solution to the associated Riemann problem,  traffic density and speed are non-negatives and bounded, the infor\-mation is not faster than traffic flow, and the shock-waves are related to breaking meanwhile rarefaction-waves are associated to acceleration. In \citep{Jiang2014,Jiang2016} the authors use the corrected pressure model, but adding on the right side of equation (\ref{ModAx:massconservation}) a non-zero function $q$, turning it into a pseudo-conservation law,  and including also in (\ref{ModAx:momentum}) the relaxation term corresponding to the desired speed $\mathbf{v}$.  This modified model was used to simulate the traffic flow and its consequences in air pollution on a two-dimensional city. Despite the pseudo-conservation of cars employed in \citep{Jiang2014,Jiang2016}, they assume that mathematical properties -arisen from the correction on pressure- demonstrated in \cite{Aw2000} still holds, being this statement a pending task for future studies.

\subsection{Porous media models and the interchange between fluid and solid phases}

A porous media could be defined in a rough manner as a domain which combines a solid part (denoted solid-phase) with a complementary part with voids filled by a fluid (denoted fluid-phase). The fluid can flow through channels connecting the voids, allowing the transport and the interchange of mass and energy between solid and fluid phases. The proportion of both phases is given by the so-called porosity, denoted here by $\epsilon$, being a dimensionless value between $0$ and $1$. Thus, in a representative elementary volume (REV) the porosity $\epsilon$ is the proportion of voids (fluid-phase), and -due the complementary nature of both phases- then ($1-\epsilon$) is the proportion of the solid phase.  It is worthwhile remarking here that the extreme cases are no more a porous media: only solid for $\epsilon=0$, or immaterial for $\epsilon=1$. 

There is a widely extensive literature in fluid dynamics on porous media, which turns the mathematical models in PDE for transport, and for the interchange of energy and mass, into a well-known matter (see \cite{das2018_book} and references therein).  However, it is important remarking here that our quantity under study in the fluid-phase is neither mass nor energy, but the fluid density itself.  This idea will be a key point when our model is formulated.  The equations proposed in this work for a compressible flow on porous media (related to continuity and momentum of traffic flow) will be formulated in detail in the following subsection. 

\subsection{A novel nonconservative traffic flow model in an urban-porous media}

\begin{figure}
\centering
\subcaptionbox{A section of the city of Guadalajara (taken from Google-Earth (2024)).\label{fig1:GDL}}{\includegraphics[width=.8\linewidth,height=.2\textheight,keepaspectratio]{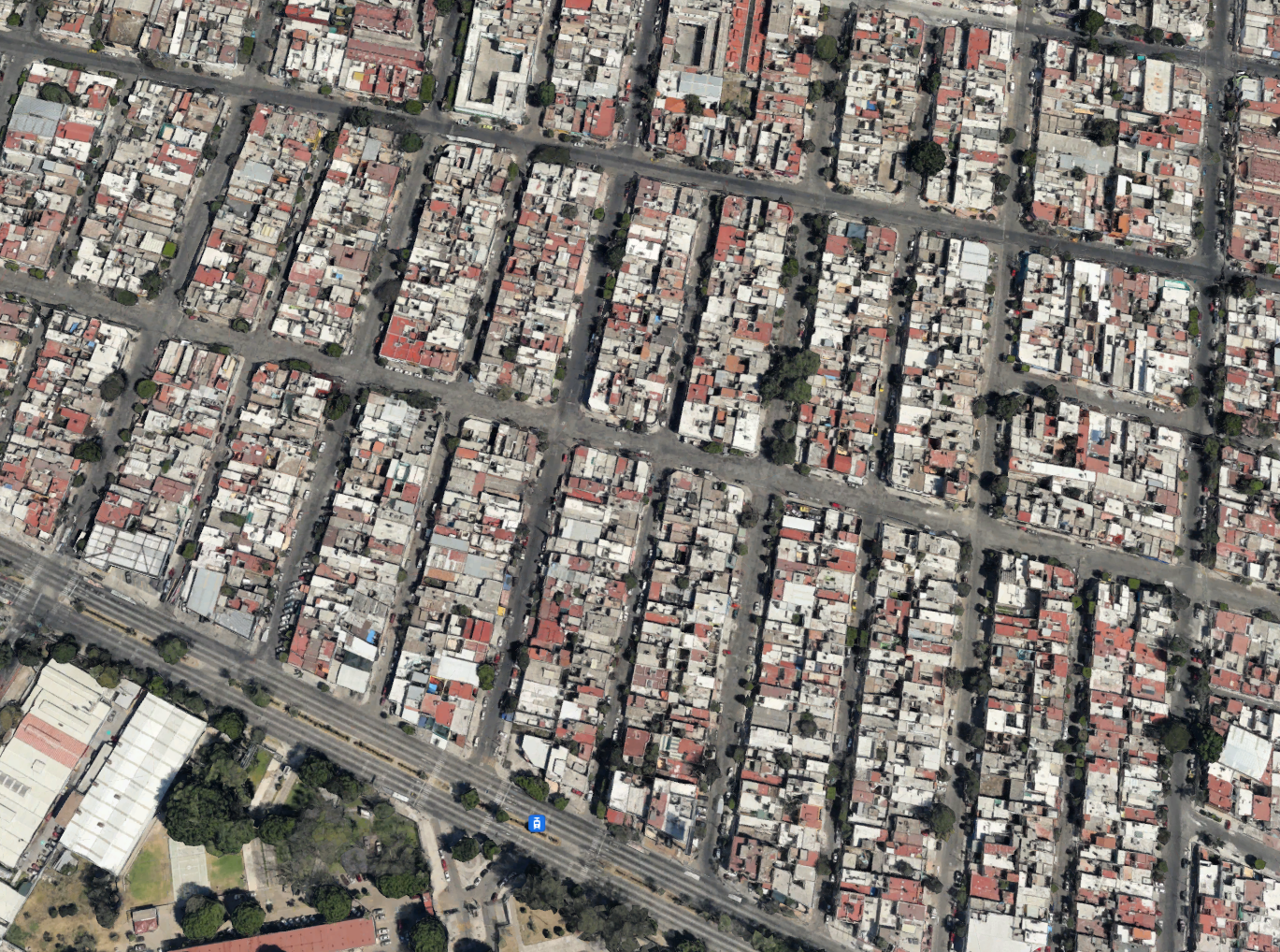}}
\subcaptionbox{Diagram of the nonconservative traffic flow model on an urban-porous city.\label{fig1:diagram}}{\includegraphics[width=.8\linewidth,height=.2\textheight,keepaspectratio]{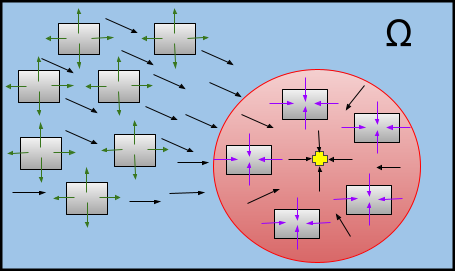}}
\caption{(\subref{fig1:GDL}) A city as an urban-porous media domain: city landscape is characterized by building-blocks delimited by streets, in this case in Guadalajara city.  (\subref{fig1:diagram}) Scheme of the porous traffic flow model: cars leave residential buildings (green arrows) on suburbia areas (in blue) to enter the city center (in red), flowing on the urban-porous media towards the attraction point represented by a yellow cross (black arrows). This attraction point belongs to a zone containing malls, working places, or schools, were cars leave the fluid phase going into the corresponding parking-spaces (magenta arrows).  \label{fig1}}
\end{figure}

Using a drone is simple to observe the city as a set of building-blocks clearly delimited by streets. This set up is similar to a porous media where the solid-phase is delimited by interconnected voids (see Figure \ref{fig1:GDL}) of its fluid-phase.  On this combination of streets and building-blocks is possible to visualize an everyday dynamic where cars leave their parking spaces, either on the street sides or inside of buildings, to flow on streets, finishing their travel for back again at a parking-space on their destination.  In order to model this dynamic, we need to specify a few assumptions on the process: 
\begin{itemize}
\item The city has its residential zones located close to city limits (suburbia), meanwhile the working-places, malls, schools, and other attraction centers are most of them in the city center. 
\item Inside the city there exist areas, such as urban parks or pedestrian zones, where the traffic flow is forbidden. 
\item Every day drivers leave their residential buildings entering the streets. Then, traffic flows inside the city through the porous medium until reaching the attraction centers, where drivers look for an off-street parking space that, for simplicity,  is assumed to have unlimited capacity.  When they find it, cars leave the streets and go inside buildings, leaving free space on streets for new drivers which continuously reach the city center (see Figure \ref{fig1:diagram}).  The case of on-street parking can be treated, for instance, with the techniques introduced in \cite{lsm2017}.
\end{itemize}

This implies a nonconservative dynamic for cars if we only follow them on the streets. However,  as already commented,  this interchange is different to that given by classical equations, because the quantity transported on the urban porous media is cars density itself.  So,  its interchange must be modeled on the continuity equation (see below equation (\ref{ModTraffic:continuidad})). Moreover, we complement this equation with two terms: a traffic demand (mass source) term depending on the building-blocks (solid-phase), and an absorption term depending on the fluid-phase (streets), on traffic density, and on an absorption rate (parking). This dynamic of mass source/sink is typical in pollution modeling \cite{sportisse2010},  and we use it to formulate the continuity equation for the traffic density.  Finally, we also include in the continuity equation small enough viscosity effects to stabilize the numerical solution.

Concerning the momentum equation,  in addition to the already introduced correction of the pressure \citep{Aw2000} to adjust density and speed to sonic velocity and the relaxation term to measure the drivers' response to a desired speed,  we also include in the model a stabilizing viscosity term.  However, although traffic flow does not reach subsonic velocities,  is actually faster than typical flows on porous media.  So, additional terms (corrections) to classical Darcy law are needed: in particular, the Brinkman correction is added to the viscosity effects, and also the so-called Forchheimer correction.  All above mentioned elements and assumptions lead us to formulate the novel noncon\-servative model for the traffic flow in the following manner.

Let $\Omega\subset\mathbb{R}^2$ be an urban-porous media with boundary segments $\Gamma = \Gamma_{w}\cup\Gamma_{bnd}$,  associated respectively to the obstacles' walls inside $\Omega$ ($\Gamma_{w}$), and to the urban limit layout ($\Gamma_{bnd}$).  Let $(0,T)$ be the time interval of interest.  Then, we look for the traffic density and velocity such that satisfy the following system of PDE:  
\begin{subequations}\label{ModTraffic}
\begin{align}
   \epsilon\dot{\rho} + \nabla\cdot (\rho\mathbf{u}) - \nabla\cdot (\epsilon\nu\nabla\rho) + \epsilon\kappa\rho = (1-\epsilon)q &\quad \mbox{ in }\Omega\times (0,T),\label{ModTraffic:continuidad}\\
   \dot{\mathbf{u}} + \frac{1}{\epsilon} (\mathbf{u}\cdot\nabla) \mathbf{u} - c^2\rho(\nabla\cdot\mathbf{u})\mathbf{1} - \frac{\mathbf{v} - \mathbf{u}}{\tau} 
   & \nonumber\\
 = \frac{\epsilon}{\rho}(\nabla\cdot\frac{\mu}{\epsilon}\nabla) \mathbf{u}  - \frac{\epsilon \mu}{\rho K} \mathbf{u} - \frac{\epsilon F}{\sqrt{K}}\Vert\mathbf{u}\Vert \mathbf{u} & \quad \mbox{ in }\Omega\times (0,T),\label{ModTraffic:momentum}   
\end{align}
\end{subequations}
where $\rho\ [{\rm veh/Km}^2]$ is the traffic density, $\mathbf{u}\ [{\rm Km/h}]$ is the local speed of traffic, $\epsilon$ is the urban porosity, $(1-\epsilon)q\ [{\rm veh/Km}^2{\rm /h}]$ is the traffic demand function associated to the building-blocks (solid-phase), $\mathbf{v}\ [{\rm Km/h}]$ is the desired speed, $c^2$ is the sonic speed, $K$ is the urban permeability, $F$ is the Forchheimer coefficient, $\nu\ [{\rm Km}^2{\rm /h}]$ is the coefficient of mass diffusion, $ \kappa\ [{\rm h}^{-1}]$ is the absorption rate (parking) from the streets (fluid-phase) to parking-spaces inside building-blocks (solid-phase), $\tau\ [{\rm h}]$ is the relaxation time, and $\mu\ [{\rm Km}^2{\rm /h}]$ is the viscosity.  

System (\ref{ModTraffic}) is complemented with suitable initial and boundary conditions. So, we impose a Neumann homogeneous condition, and a slip type condition on the two boundary segments to avoid that cars can scape by diffusion, and to enforce a tangential flow of cars when an obstacle is on their path:
\begin{subequations}\label{ModelBnd}
\begin{align}
\nabla\rho\cdot\mathbf{n} = 0 &\quad \mbox{ on } \Gamma_{w}\cup\Gamma_{bnd}, \label{ModelBnd:neumann}\\   
\mathbf{u}\cdot\mathbf{n} = 0 &\quad \mbox{ on } \Gamma_{w}, \label{ModelBnd:slip}\\
\nabla\mathbf{u} \mathbf{n} = \mathbf{0} &\quad \mbox{ on } \Gamma_{bnd}.  \label{ModelBnd:new}
\end{align}
\end{subequations}  
We must note here that above boundary condition $\nabla\mathbf{u} \mathbf{n} = \mathbf{0}$ is equivalent to conditions $\nabla u_1\cdot\mathbf{n} = \nabla u_2\cdot\mathbf{n} = 0$. 

Finally, we assume that, at the initial time $t=0$, the traffic density has a known distribution $\rho^0$ with a null traffic speed:
\begin{subequations}\label{ModIni}
\begin{align}
\rho(.,0) = \rho^0(.) & \quad \mbox{ in } \Omega \label{ModIni:densidad},\\
\mathbf{u}(.,0) = \mathbf{0} & \quad \mbox{ in } \Omega\label{ModIni:velocidad}.
\end{align}
\end{subequations}

\subsection{Modeling a desired traffic speed: a modified Eikonal equation}\label{subsec:Eikonal}

The nonconservative traffic flow model (\ref{ModTraffic}) requires fixing a desired traffic speed $\mathbf{v}(\rho)$, which could be obtained from data, drivers' patterns, or mathe\-matical models. In our case, we use the latter option as already did previously by other authors \cite{Xia2008,Jiang2015,Jiang2016}.  So, we employ an Eikonal-based approach, where the Eikonal model -giving the speed direction- is formulated in the following way:
\begin{subequations}\label{ModEikon}
\begin{align}
\Vert \nabla\phi\Vert = \frac{1}{f(\rho)} & \quad \mbox{ in }\Omega,\label{ModEikon:Eq}\\
\phi = 0 & \quad \mbox{ on } \Gamma_{exit},\label{ModEikon:bndexit}\\
\nabla\phi\cdot\mathbf{n} = 0 & \quad \mbox{ on }\Gamma_w \cup\Gamma_{bnd},\label{ModEikon:bndwall}
\end{align}
\end{subequations}
with $\phi$ representing the total instantaneous travel cost, and $f(\rho)$ the transpor\-tation cost -which depends on the traffic density. We must note that the solution $\phi$ is a potential (i.e.,  it gives us the direction $\nabla\phi$ pointing to the ``exit'' located on an inner boundary curve $\Gamma_{exit}$,  delimiting the parking zone,  where $\phi$ is null). Of course, this direction can be formulated as the normalized vector $\nabla\phi/\Vert\nabla\phi\Vert$ to be used only as a desired direction.  For instance, in \cite{Xia2008} this potential was used to drive pedestrian towards a train station directed to an exit door, meanwhile in \cite{Jiang2016} the travel cost was used to conduct the traffic flow outside the city thought a domain-hole at city center delimited by the boundary segment $\Gamma_{exit}$. 

Problem (\ref{ModEikon}) is nonlinear, but can be turn on a linear problem with the following steps (see further details in \cite{axthelm2016}). We square the equation (\ref{ModEikon:Eq}) and add a diffusive term with a suitable $\eta>0$,  so we write:
\begin{equation}
-\eta\Delta\phi + \Vert\nabla\phi\Vert^2 = \frac{1}{f^2(\rho)}.
\end{equation} 
Taking the new variable $\psi = e^{-\phi / \eta}$ (this is, $\phi = -\eta\,\ln(\psi)$), and replacing it on above equation, we can reformulate the problem with unknown $\psi$:
$$
\eta^2\psi^{-1}\Delta\psi = \frac{1}{f^2(\rho)}.
$$
Now, multiplying this equation by $\psi$, we arrive to the following linear problem, forced by $f(\rho)$ and with consistent boundary conditions:
\begin{subequations}\label{ModPsi}
\begin{align}
\eta^2\Delta\psi - \frac{1}{f^2(\rho)}\psi = 0 & \quad \mbox{ in }\Omega,\label{ModPsi:Eq}\\
\psi = 1 & \quad \mbox{ on }\Gamma_{exit},\label{ModPsi:bndexit}\\
\nabla\psi\cdot\mathbf{n} = 0 & \quad \mbox{ on }\Gamma_w\cup\Gamma_{bnd}.\label{ModPsi:bndwall}
\end{align}
\end{subequations}
This linear model (\ref{ModPsi}) maintains the idea of the ``exit'' as an inner boundary encircling the attraction area corresponding to the off-street parking zone.  However,  for our current purposes, we prefer to consider a modified model where the new Eikonal direction points to a central attraction point (corres\-ponding to the center of the parking domain). To do this, a new forcing term is added, being independent of both the density and the unknown $\psi$,  that strengthens the desired direction to target this attraction point. Thus, we rewrite the model as follows:
\begin{subequations}\label{ModPsi2}
\begin{align}
\eta^2\Delta\psi - \frac{1}{f^2(\rho)}\psi = G & \quad \mbox{ in }\Omega,\label{ModPsi2:eq}\\
\nabla\psi\cdot\mathbf{n} = 0 & \quad \mbox{ on }\Gamma_w\cup\Gamma_{bnd},\label{ModPsi2:bndwall}
\end{align}
\end{subequations}
where the second member term $G$ is chosen in such a way that $\psi$ tends to 1 as $\mathbf{x}$ tends to the attraction point (neglecting in this way the boundary condition (\ref{ModPsi:bndexit}) on $\Gamma_{exit}$). Obviously, in this way, $\phi$ tends to 0 as $\mathbf{x}$ tends to the attraction point (substituting condition (\ref{ModEikon:bndexit})).
Further details on $G$ can be found in below section \ref{se4}.

Therefore, to obtain the desired traffic speed $\mathbf{v}(\rho)$ we solve first the linear system (\ref{ModPsi2}), then we make $\phi = -\eta\,\ln(\psi) $,  we define the normalized desired direction $\nabla\phi/\Vert\nabla\phi\Vert$ and, finally, we include the travel cost $f(\rho)$ to formulate the desired traffic speed:
\begin{equation}\label{Eq:DesiredSpeed}
\mathbf{v}(\rho) = -f(\rho)\frac{\nabla\phi}{\Vert\nabla\phi\Vert}.
\end{equation}

\section{Numerical solution}\label{sec:numsol}

To solve our model,  we take advantage of Brinkman correction (viscosity) and density diffusion presence to implement a numerical scheme consisting of a combination of standard finite elements of Lagrange $P_1$ type in space with an explicit time marching scheme.  Among the numerous explicit time methods,  we use the so-called strong stability-preserving method with the aim of avoiding instabilities on the numerical solution.  Both strategies are shown in below subsections,  beginning with the variational formulation and the space discretization.     

\subsection{Variational formulation and the semi-discretized problem}

In order to obtain a numerical solution of the problem,  we need first to address the variational form of our model.  So,  multiplying equations (\ref{ModTraffic:continuidad}),  (\ref{ModPsi2:eq}) and (\ref{ModTraffic:momentum}) by respective test functions $(v,v,\mathbf{w})\in V\times V\times \mathbf{W}$, for suitable functional spaces $V$ and $\mathbf{W}$,  and integrating on $\Omega$ we arrive to the variational formulation:
\begin{align*}
\int_{\Omega} \epsilon\dot{\rho} v + \int_{\Omega} \nabla\cdot(\rho\mathbf{u}) v - \int_{\Omega} \nabla\cdot(\epsilon\nu\nabla\rho) v + \int_{\Omega} \epsilon\kappa\rho v = \int_{\Omega} (1-\epsilon)q v, \\
\int_\Omega \eta^2 \Delta\psi v - \int_\Omega \frac{1}{f^2(\rho)} \psi v = \int_\Omega Gv, \\
\int_{\Omega} \dot{\mathbf{u}}\odot \mathbf{w} + \int_{\Omega} \frac{1}{\epsilon}(\mathbf{u}\cdot\nabla) \mathbf{u} \odot \mathbf{w} - \int_{\Omega} c^2\rho(\nabla\cdot\mathbf{u)\mathbf{1}}\odot \mathbf{w} - \int_{\Omega}\frac{\mathbf{v} - \mathbf{u}}{\tau}\odot \mathbf{w} \\
 = \int_{\Omega} \frac{\epsilon}{\rho} (\nabla\cdot\frac{\mu}{\epsilon}\nabla) \mathbf{u}\odot \mathbf{w} 
   - \int_{\Omega}\frac{\epsilon\mu}{\rho K} \mathbf{u}\odot \mathbf{w} - \int_{\Omega} \frac{\epsilon F}{\sqrt{K}}\Vert\mathbf{u}\Vert \mathbf{u}\odot \mathbf{w} ,
\end{align*} 
where the products $\mathbf{u\odot w}$ represent the component-wise multiplication (also known as Hadamard product). 

In accordance with the standard finite element method, the above system suggests following integration by parts to apply boundary conditions: 
\begin{align*}
-\int_{\Omega} \nabla\cdot(\epsilon\mu\nabla\rho) v = \int_{\Omega} \epsilon\nu\nabla\rho\cdot\nabla v - \int_{\Gamma_w\cup\Gamma_{bnd}} \epsilon\nu(\nabla\rho\cdot\mathbf{n})v,\\
\int_\Omega \eta^2 \Delta\psi v = -\int_\Omega \eta^2\nabla\psi\cdot\nabla v + \int_{\Gamma_{w}\bigcup\Gamma_{bnd}} \eta^2\nabla\psi\cdot\mathbf{n},\\
\int_{\Omega}\frac{\epsilon}{\rho}(\nabla\cdot\frac{\mu}{\epsilon}\nabla) \mathbf{u} \odot \mathbf{w} = -\int_{\Omega}\frac{\mu}{\epsilon}\nabla \mathbf{u} \odot\nabla(\frac{\epsilon}{\rho}\mathbf{w}) + 
\int_{\Gamma_w\cup\Gamma_{bnd}} \frac{\mu}{\rho}(\nabla\mathbf{u}\mathbf{n}) \odot \mathbf{w}.
\end{align*}
Applying now conditions (\ref{ModelBnd:neumann}), (\ref{ModelBnd:new}) and (\ref{ModPsi2:bndwall}) allows us to drop the integral on $\Gamma_w\cup\Gamma_{bnd}$ in the first equation (ensuring cars do not leave the city by diffusion) and in the second equation.  However,  we are able only to drop the integral on the city limit $\Gamma_{bnd}$ in the third equation. 

On the other hand,  the slip condition (\ref{ModelBnd:slip}) requires a different type of arguing.  Thus, integrating by parts the terms related to mass convection and pressure correction (assuming $c^2$ as constant) we have:
\begin{align*}
\int_{\Omega}\nabla\cdot(\rho\mathbf{u}) v = -\int_{\Omega} \rho(\mathbf{u} \cdot\nabla v) + \int_{\Gamma_w\cup\Gamma_{bnd}} \rho(\mathbf{u}\cdot\mathbf{n})v,\\
- c^2 \int_{\Omega}\rho(\nabla\cdot\mathbf{u})\mathbf{1}\odot \mathbf{w} = c^2\int_{\Omega} \nabla(\rho \mathbf{w}) \mathbf{u} - c^2\int_{\Gamma_w\cup\Gamma_{bnd}} \rho(\mathbf{u}\cdot\mathbf{n})\mathbf{1}\odot\mathbf{w},
\end{align*}
allowing us to apply the slip condition on $\Gamma_w$, remaining only boundary integrals on $\Gamma_{bnd}$. 

So, we can formulate the follow variational problem: Find $(\rho,\psi,\mathbf{u})\in V\times V\times \mathbf{W}$ such that for all $(v,v,\mathbf{w})\in V\times V\times \mathbf{W}$ satisfies the system
\begin{subequations}
\begin{align*}
\int_{\Omega}\epsilon\dot{\rho} v -\int_{\Omega} \rho(\mathbf{u}\cdot\nabla v_h) + \int_{\Gamma_{bnd}} \rho(\mathbf{u}\cdot\mathbf{n})v + \int_{\Omega} \epsilon\nu\nabla\rho\cdot\nabla v  \\ 
+ \int_{\Omega}\epsilon\kappa\rho v = \int_{\Omega}(1-\epsilon)q v,\\
\int_\Omega \eta^2\nabla\psi\cdot\nabla v - \int_\Omega \frac{1}{f^2(\rho)} \psi v = \int_\Omega Gv,\\
\int_{\Omega} \dot{\mathbf{u}}\odot\mathbf{w} + \int_{\Omega}\frac{1}{\epsilon}(\mathbf{u}\cdot\nabla)\mathbf{u}\odot\mathbf{w} + c^2\int_{\Omega}  \nabla(\rho \mathbf{w})\mathbf{u} \nonumber \\ 
- c^2\int_{\Gamma_{bnd}} \rho(\mathbf{u}\cdot\mathbf{n})\mathbf{1}\odot\mathbf{w} - \int_{\Omega} \frac{\mathbf{v} - \mathbf{u}}{\tau}\odot\mathbf{w}  = - \int_{\Omega}\frac{\mu}{\epsilon}\nabla \mathbf{u} \odot\nabla(\frac{\epsilon}{\rho}\mathbf{w})  \\
+ \int_{\Gamma_w}\frac{\mu}{\epsilon}\nabla\mathbf{un}\odot\mathbf{w} -\int_{\Omega} \frac{\epsilon\mu}{\rho K}\mathbf{u}\odot\mathbf{w} - \int_{\Omega}\frac{\epsilon F}{\sqrt{K}}\Vert\mathbf{u}\Vert \mathbf{u}\odot\mathbf{w}.
\end{align*}
\end{subequations}
Now, let $\Omega_h$ be a polygonal approximation of $\Omega$, and let $\mathcal{K}_h$ be a suitable triangular mesh of $\Omega_h$ with $n_t$ triangular elements ${K}_i,\, i=1,2,\ldots,n_t,$ and $n_v$ nodes $\mathbf{x}_v, v=1,\ldots,n_p$.  So, let us consider the finite-dimensional spaces $V_h$ and $\mathbf{W}_h$ associated to the mesh $\tau_h$, corresponding to functional spaces $V$ and $\mathbf{W}$,  defined as follows: 
\begin{align*}
V_h = \{ v\in \mathcal{C}(\bar\Omega): v\vert_{K_i}\in P_1, \forall i = 1, \dots, n_t, \nabla v\cdot\mathbf{n} = 0 \mbox{ on } \Gamma_{w}\cup\Gamma_{bnd}\},\\
\mathbf{W}_h = \left\{\mathbf{w}\in [\mathcal{C}(\bar\Omega)]^2: \mathbf{w}\vert_{K_i}\in [P_1]^2, \forall i = 1, \dots, n_t, \nabla\mathbf{w}\mathbf{n}=\mathbf{0} \mbox{ on } \Gamma_{bnd},\right. \\ 
\left. \mbox{ and } \mathbf{w}\cdot\mathbf{n}=0 \mbox{ on } \Gamma_{w} \right\}.
\end{align*}
Taking $\rho^0_h \in V_h$ such that $\rho^0_h\approx\rho^0$,  we can formulate the semi-discrete problem: 
Given the initial conditions $\rho_h(.,0)=\rho^0_h, \mathbf{u}_h(.,0)=\mathbf{0}$, find the finite element functions $(\rho_h,\psi_h,\mathbf{u}_h)\in V_h\times V_h\times \mathbf{W}_h$ such that for all $(v_h,v_h,\mathbf{w}_h)\in V_h\times V_h\times \mathbf{W}_h$, satisfies the system:
\begin{subequations}\label{ProbSemidiscreto} \nonumber
\begin{align}
\int_{\Omega}\epsilon\dot{\rho}_h v_h -\int_{\Omega} \rho_h(\mathbf{u}_h\cdot\nabla v_h) + \int_{\Gamma_{bnd}} \rho_h(\mathbf{u}_h\cdot\mathbf{n})v_h + \int_{\Omega} \epsilon\nu\nabla\rho_h\cdot\nabla v_h  \nonumber\\
+ \int_{\Omega}\epsilon\kappa\rho_h v_h = \int_{\Omega}(1-\epsilon)q v_h,\\
\int_\Omega \eta^2\nabla\psi_h\cdot\nabla v_h - \int_\Omega \frac{1}{f^2(\rho_h)} \psi_h v_h = \int_\Omega G_hv_h,\\
\int_{\Omega} \dot{\mathbf{u}}_h\odot\mathbf{w}_h + \int_{\Omega}\frac{1}{\epsilon}(\mathbf{u}_h\cdot\nabla)\mathbf{u}_h\odot\mathbf{w}_h + c^2\int_{\Omega}\nabla(\rho_h\mathbf{w}_h)\mathbf{u}_h  \nonumber\\ 
- c^2\int_{\Gamma_{bnd}} \rho(\mathbf{u}_h\cdot\mathbf{n})\mathbf{1}\odot\mathbf{w}_h - \int_{\Omega} \frac{\mathbf{v}_h - \mathbf{u}_h}{\tau}\odot\mathbf{w}_h  = -\int_{\Omega}\frac{\mu}{\epsilon}\nabla \mathbf{u}_h \odot\nabla(\frac{\epsilon}{\rho_h}\mathbf{w}_h)  \nonumber\\
+ \int_{\Gamma_w}\frac{\mu}{\epsilon}\nabla\mathbf{u}_h\mathbf{n}\odot\mathbf{w}_h - \int_{\Omega} \frac{\epsilon\mu}{\rho_h K}\mathbf{u}_h\odot\mathbf{w}_h - \int_{\Omega}\frac{\epsilon F}{\sqrt{K}}\Vert\mathbf{u}_h\Vert \mathbf{u}_h\odot\mathbf{w}_h.
\end{align}
\end{subequations}
Then,  following the standard finite element process \cite{Knaber2003,Ern2021},  we consider the basis $\mathcal{B}_{V_h} = \{\tilde{v}_i\}^{n_p}_{i=1}$ and $\mathcal{B}_{\mathbf{W}_h} = \{\tilde{\mathbf{w}}_i\}^{n_p}_{i=1}$ of the subspaces $V_h$ and $\mathbf{W}_h$ respectively, being its elements such that $\tilde{v}_i(\mathbf{x}_j) = \delta^j_i$ and $\tilde{\mathbf{w}}_i(\mathbf{x}_j) = [\delta^j_i, \delta^j_i]'$ for all $i,j \in \{1,\ldots,n_v\}$, and given that $(\rho_h,\psi_h,\mathbf{w}_h)$ belongs to $V_h\times V_h\times\mathbf{W}_h$ then admits the following expansions: 
\begin{equation}
\rho_h = \sum^{n_p}_{i=1} \xi_i(t)\tilde{v}_i(\mathbf{x}),\quad \psi_h = \sum^{n_p}_{i=1}\zeta_i\tilde{v}_i(\mathbf{x}),\quad \mathbf{u}_h=\sum^{n_p}_{i=1}\boldsymbol\upsilon_i(t)\odot\tilde{\mathbf{w}}_i(\mathbf{x}),
\end{equation} 
where it must be noted that the expansion coefficients are $\xi_i(t) = \rho_h(\mathbf{x}_i,t)$,  $\zeta_i = \psi_h(\mathbf{x}_i)$ and $\boldsymbol{\upsilon}_i(t) = \mathbf{u}_h(\mathbf{x}_i,t)$. Now we take as test functions the elements of the basis $v_h = \tilde{v}_j, \mathbf{w}_h = \tilde{\mathbf{w}}_j$ for each index $j\in\{1,\ldots,n_v\}$.  So, substituting the unknowns' expansions and the chosen test functions in the variational problem it turns into the ordinary differential equations (ODE) system:    
\begin{subequations}\label{ODEsys}
\begin{align}
\sum^{n_v}_{i=1}\dot{\xi}_i(t)\int_{\Omega} \epsilon\tilde{v}_i\tilde{v}_j + \sum^{n_v}_{i=1}\xi_i(t)[-\int_\Omega \tilde{v}_i(\mathbf{u}_h\cdot\nabla\tilde{v}_j) + \int_{\Gamma_{bnd}} \tilde{v}_i(\mathbf{u}_h\cdot\mathbf{n})\tilde{v}_j  \nonumber\\
+ \int_\Omega \epsilon\nu\nabla\tilde{v}_i\cdot\nabla\tilde{v}_j + \int_\Omega \epsilon\kappa\tilde{v}_i\tilde{v}_j] = \int_\Omega (1-\epsilon)q\tilde{v}_j , \label{ODEsys:densidad}\\
\sum^{n_v}_{i=1}\zeta_i [\int_\Omega \eta^2\nabla\tilde{v}_i\cdot\nabla\tilde{v}_j - \int_\Omega \frac{1}{f^2(\rho_h)} \tilde{v}_i\tilde{v}_j] = \int_\Omega G\tilde{v}_j, \label{ODEsys:eikonal}\\ 
\sum^{n_v}_{i=1} \dot{\boldsymbol{\upsilon}}_i(t) \int_\Omega \tilde{\mathbf{w}}_i\odot\tilde{\mathbf{w}}_j + \sum^{n_v}_{i=1} \boldsymbol{\upsilon}_i(t)[\int_\Omega \frac{1}{\epsilon}(\mathbf{u}_h\cdot\nabla)\tilde{\mathbf{w}}_i\odot\tilde{\mathbf{w}}_j  \nonumber\\ 
+ \int_\Omega\frac{1}{\tau}\tilde{\mathbf{w}}_i\odot\tilde{\mathbf{w}}_j + \int_\Omega\frac{\mu}{\epsilon}\nabla\tilde{\mathbf{w}}_i\odot\nabla\frac{\epsilon}{\rho_h}\tilde{\mathbf{w}}_j - \int_{\Gamma_{bnd}} \frac{\mu}{\epsilon}\nabla\tilde{\mathbf{w}}_i\mathbf{n}\odot\tilde{\mathbf{w}}_j  \nonumber\\ 
+ \int_\Omega \frac{\epsilon\mu}{\rho_h K} \tilde{\mathbf{w}}_i\odot\tilde{\mathbf{w}}_j + \int_\Omega \epsilon\frac{F}{\sqrt{K}} \Vert\mathbf{u}_h\Vert\tilde{\mathbf{w}}_i\odot\tilde{\mathbf{w}}_j] =  - c^2\int_\Omega \nabla(\rho_h\tilde{\mathbf{w}}_j)\mathbf{u}_h  \nonumber\\
+ c^2\int_{\Gamma_{bnd}} \rho_h(\mathbf{u}_h\cdot\mathbf{n}) \mathbf{1}\odot\tilde{\mathbf{w}}_j  + \int_\Omega \frac{1}{\tau}\mathbf{v}(\rho_h)\odot\tilde{\mathbf{w}}_j. \label{ODEsys:velocidad}
\end{align}
\end{subequations}
Notice that we conserve the unknowns $\mathbf{u}_h$ and $\rho_h$ to be treated as a given data to choose a suitable marching time scheme of explicit type. Also, the system presents two integrals that involve the gradient of functions' product, $\nabla(\frac{\epsilon}{\rho_h}\tilde{\mathbf{w}}_j)$ and $\nabla(\rho_h\tilde{\mathbf{w}}_j)$,  making necessary to simplify them to reduce those integrals to well-known integrals of the finite element method. This can be done on the discrete level using the mean values $\overline{\rho}_h\vert_{K_i}$ and $\bar{\epsilon}\vert_{K_i}$ over each triangle $K_i$
\begin{align*}
\int_{\Omega} \frac{\mu}{\epsilon}\nabla \mathbf{w}_i \odot\nabla(\frac{\epsilon}{\rho_h}\mathbf{w}_j) \approx \sum_{i=1}^{n_t}\frac{1}{\overline{\rho}_h\vert_{K_i}}\int_{K_i} \mu\nabla \mathbf{w}_i \odot\nabla \mathbf{w}_j,\\
c^2\int_{\Omega} \nabla(\rho_h \mathbf{w}_j)\mathbf{u}_h \approx \sum^{n_t}_{i=1}c^2\overline{\rho}_h\vert_{K_i}\int_{K_i}  \nabla \mathbf{w}_j\mathbf{u}_h.
\end{align*}
It is well-known that the assembling process of the finite element method turn the integrals in matrices and vectors, for example the called mass-matrix $\mathsf{M}_h$ from equation (\ref{ODEsys:densidad}) defined by 
$
(\mathsf{M}_h)_{i,j} = \int_\Omega \tilde{v}_i\tilde{v}_j
$
and the demand vector $\mathsf{b_h}$ given by $(\mathsf{b_h})_j = \int_\Omega (1-\epsilon)q\tilde{v}_j$, with dimensions $n_v\times n_v$ and $n_v\times 1$ respectively. However, we must note than in equation (\ref{ODEsys:velocidad}) the assembled matrices and vectors have dimensions $2n_v\times 2n_v$ and $2n_v\times 1$.  Indeed, the assembling process on integrals like $\int_\Omega \tilde{\mathbf{w}}_i\odot\tilde{\mathbf{w}}_j$ and $\int_\Omega \frac{1}{\tau}\mathbf{v}(\rho_h)\odot\tilde{\mathbf{w}}_j$ gives us a block matrix and a vector:
$$
\tilde{\mathsf{M}}_h = \begin{bmatrix}
                       \mathbf{M} & 0 \\
                       0 & \mathbf{M}    
                       \end{bmatrix},\quad \tilde{\mathsf{c}}_h = \begin{bmatrix}
                                                                   \mathbf{c}_1\\
                                                                   \mathbf{c}_2                                                 
                                                                  \end{bmatrix},
$$
where $\mathbf{M}$ is the standard mass-matrix obtained after assembling the integral $(\mathbf{M})_{i,j}=\int_\Omega \tilde{w}_i\tilde{w}_j$ of basis vector function $\tilde{\mathbf{w}}_j$ components, and $(\mathbf{c}_k)_j = \int_\Omega \frac{1}{\tau} v_k(\rho_h)\tilde{w}_j$ for $k=1,2$. For the sake of simplicity we omit here the list of all matrices and vectors, but each one of them could be easily identified by comparing the ODE system (\ref{ODEsys}) with its matrix formulation (\ref{ODEsysMatrix}). 
 
Taking now the unknowns as the vectorial functions $\boldsymbol{\xi}(t) = [\xi_1(t),\ldots,\xi_{n_v}(t)]$, $ \boldsymbol{\upsilon}(t) = [\boldsymbol{\upsilon}_1(t),\ldots,\boldsymbol{\upsilon}_{n_v}(t)]$ and $\boldsymbol{\zeta} = [\zeta_1,\ldots,\zeta_{n_v}]$ along with the matrices and vectors from the ODE system (\ref{ODEsys}),  we can formulate the initial value problem in its matrix form:

Given the initial states $\boldsymbol{\xi}(0) = [(\rho^0_h)_1,\ldots,(\rho^0_h)_{n_v}]'$ and $\boldsymbol{\upsilon}(0) = [\mathbf{0},\ldots,\mathbf{0}]'$, find the vectors $\boldsymbol{\xi}(t), \boldsymbol{\zeta}, \boldsymbol{\upsilon}(t)$ such that satisfy the following system:
\begin{subequations}\label{ODEsysMatrix}
\begin{align}
\mathsf{M}_h \dot{\boldsymbol{\xi}}(t) + [-\mathsf{C}_h(\boldsymbol{\upsilon}) + \mathsf{B}_h(\boldsymbol{\upsilon}) + \mathsf{R}_h + \mathsf{K}_h]\boldsymbol{\xi}(t) = \mathsf{b}_h(q),\label{ODEsysMatrix:densidad}\\
[\hat{\mathsf{R}}_h - \hat{\mathsf{M}}_h(\boldsymbol{\xi})]\boldsymbol{\zeta} = \hat{\mathsf{b}}_h(G),\label{ODEsysMatrix:eikonal}\\
\tilde{\mathsf{M}}_h\dot{\boldsymbol{\upsilon}}(t) + [\tilde{\mathsf{C}}_h(\boldsymbol{\upsilon}) + \tilde{\mathsf{N}}_h + \tilde{\mathsf{R}}_h - \tilde{\mathsf{B}}_h + \tilde{\mathsf{D}}_h + \tilde{\mathsf{F}}_h(\boldsymbol{\upsilon})]\boldsymbol{\upsilon}(t)  \nonumber\\
= - c^2\tilde{\mathsf{a}}_h(\boldsymbol{\xi},\boldsymbol{\upsilon})+c^2\tilde{\mathsf{b}}_h(\boldsymbol{\xi},\boldsymbol{\upsilon}) + \tilde{\mathsf{c}}_h(\mathbf{v}(\boldsymbol{\xi}; \boldsymbol{\zeta})).\label{ODEsysMatrix:velocity}
\end{align}
\end{subequations}
As previously commented,  we have here the mass-matrices, the stiffness-matrices and the convective-matrices, but also the matrices and vectors corresponding to ensembles on the boundary and those corresponding to the Darcy law corrections. 
\subsection{Time discretization with a Strong-Stability Preserving marching time scheme}
By introducing the nodal basis of the finite-dimension subspaces $V_h$ and $W_h$ we have obtained a system of ODE with time as independent variable.  This system can be formulated in a compact form $\dot{\mathbf{U}} = \mathcal{L}(\mathbf{U})$, for the unknowns' vector $\mathbf{U} = (\boldsymbol{\xi},\boldsymbol{\upsilon})$ and the differential operator with components $\mathcal{L} = (\mathcal{F},\mathcal{G})$. These operators $(\mathcal{F},\mathcal{G})$ can be defined in a straightforward manner by solving the system (\ref{ODEsysMatrix}) for $\dot{\boldsymbol{\xi}}$ and $\dot{\boldsymbol{\upsilon}}$, respectively.   Therefore,  we have:
\begin{subequations}
\begin{align}
\mathcal{F}(\boldsymbol{\xi},\boldsymbol{\upsilon})=\mathsf{M}^{-1}_h\{[\mathsf{C}_h(\boldsymbol{\upsilon}) -\mathsf{B}_h(\boldsymbol{\upsilon}) - \mathsf{R}_h - \mathsf{K}_h]\boldsymbol{\xi} + \mathsf{b}_h(q)\}, \\
\mathcal{G}(\boldsymbol{\xi},\boldsymbol{\upsilon};\boldsymbol{\zeta}) = \tilde{\mathsf{M}}^{-1}_h\{[-\tilde{\mathsf{C}}_h(\boldsymbol{\upsilon}) - \tilde{\mathsf{N}}_h - \tilde{\mathsf{R}}_h + \tilde{\mathsf{B}}_h - \tilde{\mathsf{D}}_h - \tilde{\mathsf{F}}_h(\boldsymbol{\upsilon})]\boldsymbol{\upsilon}(t) \nonumber\\ 
-c^2\tilde{\mathsf{a}}_h(\boldsymbol{\xi},\boldsymbol{\upsilon})+c^2\tilde{\mathsf{b}}_h(\boldsymbol{\xi},\boldsymbol{\upsilon}) + \tilde{\mathsf{c}}_h(\mathbf{v}(\boldsymbol{\xi};\boldsymbol{\zeta}))\}.
\end{align} 
\end{subequations}

It is worthwhile remarking here that Eikonal equation does not present a time derivative and could be solved directly as: 
\begin{equation}
\boldsymbol{\zeta} = [\hat{\mathsf{R}}_h - \hat{\mathsf{M}}_h(\boldsymbol{\xi})]^{-1}\hat{\mathsf{b}}_h(G)
\end{equation} 
for any given vectorial function $\boldsymbol{\xi}$.

Let us introduce now a standard time discretization: Taking $N\in\mathbb{Z}^+$,  we define the time step $\delta_t = T/N$, and fix the following time instants $t^n = n\delta_t$, for $n=0,\ldots,N$. Also, we denote as $\boldsymbol{\xi}^n$ and $\boldsymbol{\upsilon}^n$ the evaluations $\boldsymbol{\xi}(t^n)$ and $\boldsymbol{\upsilon}(t^n)$,  respectively.  

Now, in agree with \cite{gottlieb_ketcheson2009}, we apply a forward time marching scheme of strong-stability preserving (SSP) type (a particular type of predictor-corrector scheme), with the aim of preventing any possible instabilities:
Given the initial states of traffic density and speed, $\boldsymbol{\xi}^0 =  [(\rho^0_h)_1,\ldots,(\rho^0_h)_{n_v}]'$,  $\boldsymbol{\upsilon}^0 = [\mathbf{0},\ldots,\mathbf{0}]'$, we look for the predictor $(\boldsymbol{\xi}^*,\boldsymbol{\upsilon}^*)$ such that satisfies the following forward scheme: 
\begin{subequations}\label{SSP-predict}
\begin{align}
\boldsymbol{\zeta} = [\hat{\mathsf{R}}_h - \hat{\mathsf{M}}_h(\boldsymbol{\xi}^n)]^{-1}\hat{\mathsf{b}}_h(G),\\
\boldsymbol{\xi}^* = \boldsymbol{\xi}^n - \delta_t \, \mathcal{F}(\boldsymbol{\xi}^n,\boldsymbol{\upsilon}^n),\\
\boldsymbol{\upsilon}^* = \boldsymbol{\upsilon}^n - \delta_t \, \mathcal{G}(\boldsymbol{\xi}^n,\boldsymbol{\upsilon}^n;\boldsymbol{\zeta}).
\end{align}
\end{subequations}
Once computed  $(\boldsymbol{\xi}^*,\boldsymbol{\upsilon}^*)$ by (\ref{SSP-predict}), we update $(\boldsymbol{\xi}^{n+1},\boldsymbol{\upsilon}^{n+1})$ with the following optimal coefficient scheme: 
\begin{subequations}\label{SSP-update}
\begin{align}
\boldsymbol{\zeta}^* = [\hat{\mathsf{R}}_h - \hat{\mathsf{M}}_h(\boldsymbol{\xi}^*)]^{-1}\hat{\mathsf{b}}_h(G),\\
\boldsymbol{\xi}^{n+1} = \boldsymbol{\xi}^n + \frac{1}{2}\boldsymbol{\xi}^* + \frac{\delta_t}{2} \, \mathcal{F}(\boldsymbol{\xi}^*,\boldsymbol{\upsilon}^*),\\
\boldsymbol{\upsilon}^{n+1} = \boldsymbol{\upsilon}^n + \frac{1}{2}\boldsymbol{\upsilon}^* + \frac{\delta_t}{2} \, \mathcal{G}(\boldsymbol{\xi}^*,\boldsymbol{\upsilon}^*; \boldsymbol{\zeta}^*),
\end{align}
\end{subequations}
for each $n=0,1,\ldots,N-1$.
Solving the SSP scheme and getting $\{\boldsymbol{\xi}^n\}^N_{n=0}$ and $\{\boldsymbol{\upsilon}^n\}^N_{n=0}$ is equivalent to compute the vehicular density and its speed for each node and instant time, indeed, we can formulate the numerical solutions of the semi-discrete problem 
as follows: $\rho_h(\mathbf{x}_i,t^n) = \boldsymbol{\xi}^n_i$, $\psi_h(\mathbf{x}_i) = \boldsymbol{\zeta}_i(\boldsymbol{\xi}^n)$ and $\mathbf{u}_h(\mathbf{x}_i,t^n) = \boldsymbol{\upsilon}^n_i$ for any $i=1,\ldots,n_v$,  and for each instant $n=0,\ldots,N$.

\section{Numerical experiments} \label{se4}

\subsection{The urban porous domain: the metropolitan zone of Guadalajara}

\begin{figure}
\centering
\subcaptionbox{The polygonal domain $\Omega$ of the part of the GMA under study.  \label{fig2:ZMGDL}}{\includegraphics[width=.6\linewidth,height=.2\textheight,keepaspectratio]{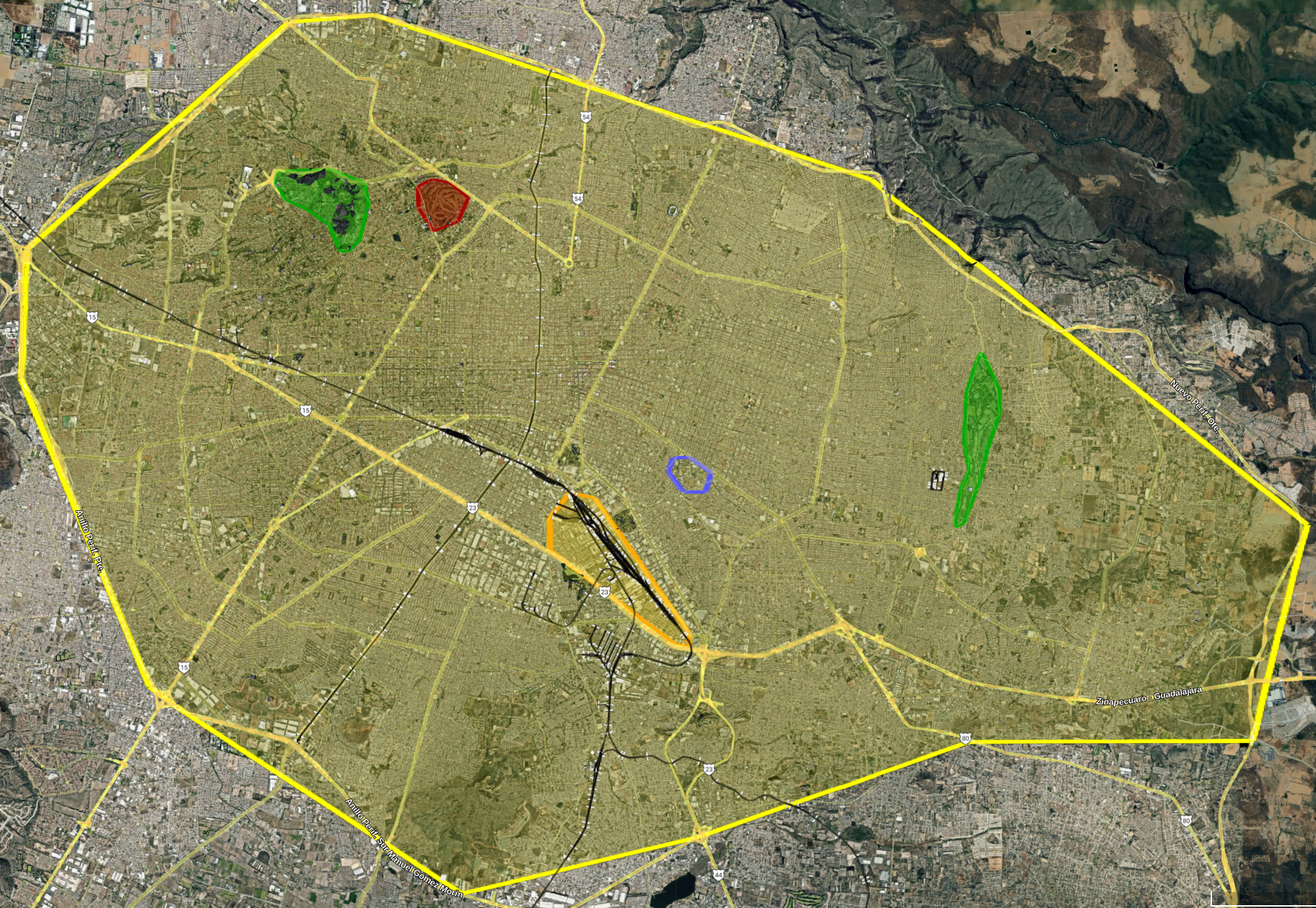}}%
\hspace{0.1cm}
\subcaptionbox{Triangular mesh of the domain $\Omega_h=\Omega$.   \label{fig2:mesh}}{\includegraphics[width=.8\linewidth,height=.22\textheight,keepaspectratio]{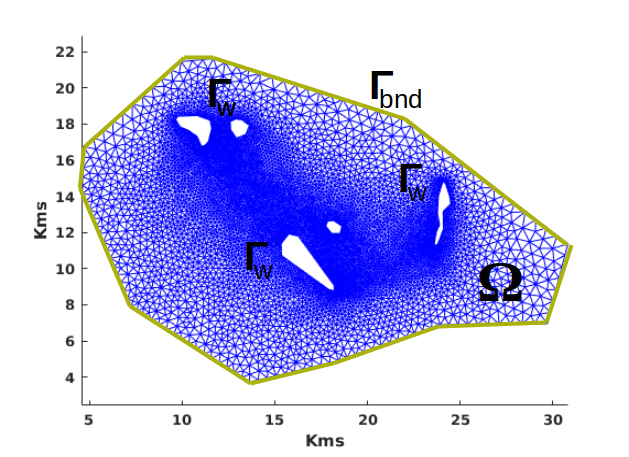}}
\caption{(\subref{fig2:ZMGDL}) The numerical domain is inspired by the metropolitan zone of Guadalajara (Mexico).  The satellite image of the zone under study shows the urban limit in yellow, and the different obstacles in green, red,  purple,  and orange (Google Earth, 2024).   (\subref{fig2:mesh}) Triangular mesh $\mathcal{K}_h$ used for the domain, depicting the urban limit boundary $\Gamma_{bnd}$,  and the obstacle walls $\Gamma_w$.}\label{fig2}
\end{figure}

The domain in our numerical tests is inspired by the Guadalajara Metro\-politan Area (GMA), in Mexico, with 5 millions of inhabitants, over 2 millions of cars of all kinds,  and only 3 metro lines.  Due to the highly irregular borders of the full metropolis, we center our study in the part of Guadalajara city that remains inside the circumvallation road surrounding its central part,  and for simplicity only a small number of obstacles inside the domain are considered: two parks, a golf club, a university campus, and an industrial zone (see Figure \ref{fig2:ZMGDL}). 

This polygonal approximation of the city $\Omega$ has been meshed employing the package GMSH \cite{gmsh2009},  using 15357 triangular elements with 7863 nodes (see Figure \ref{fig2:mesh}), with smaller triangular elements close to obstacles.  We must recall here that the boundary is divided into two parts: the city limit layout $\Gamma_{bnd}$, and the obstacle walls $\Gamma_w$ (as shown in Figure \ref{fig2:mesh}).  For the area covered by buildings and streets it was assumed a concentric structure, so that the center presents more area covered by buildings and less by streets.  In contrast, the buildings concentration decays conforming the urban limit is reached,  where more area is covered by streets and less by buildings. To model this urban configuration,  we use a simple Gaussian distribution adjusted to a scale of accepted porosity values for a city \cite{hu2012} (that is, the porosity $\epsilon$ will belong to the admissibility interval $[0.38,0.82]$) (see two different examples in Figures \ref{fig4:CiudadDensa} and \ref{fig4:CiudadDispersa}). Moreover,  we must also remark that, despite the high building density in the city center, this corresponds mainly to offices, stores, and working places.  Thus, the bulk of cars should come from the residential zones in urban environs (suburbia).

In the numerical experiments presented here, we use the following values of the model parameters (unless a new value is specified): maximum speed $U_{max}=50 \ {\rm Km/h}$,  maximum density $\rho_{max}=2000\ {\rm veh/Km}^2$, traffic demand $q = 0$, interchange rate from streets to parking zone $\kappa_{max} = 18.0 \ {\rm h}^{-1}$, density diffusion $\nu = 1.25\ {\rm Km}^2{\rm /h}$, viscosity $\mu = 3.6\times 10^{-8} \ {\rm Km}^2{\rm /h}$ (for stabilization propose), and relaxation time $\tau = 0.009 \ {\rm h}$. The simulation time was $T = 0.5\ {\rm h}$, with a time step length $\delta_t = 0.5\times 10^{-3} \ {\rm h}$.  Finally, we define the transportation cost $f(\rho)$ as the following relation: 
\begin{equation}
f(\rho) = U_{max}\big(1-\frac{\rho}{\rho_{max}}\big),
\end{equation}
which has been previously used for traffic and pedestrian flow models \cite{Xia2008,treiber2013,Jiang2016}, guaranteeing that traffic density does not exceed $\rho_{max}$.

\begin{figure}
\centering
\subcaptionbox{Attraction forcing function $G$ \label{fig3:G}}{\includegraphics[width=.8\linewidth,height=.2\textheight,keepaspectratio]{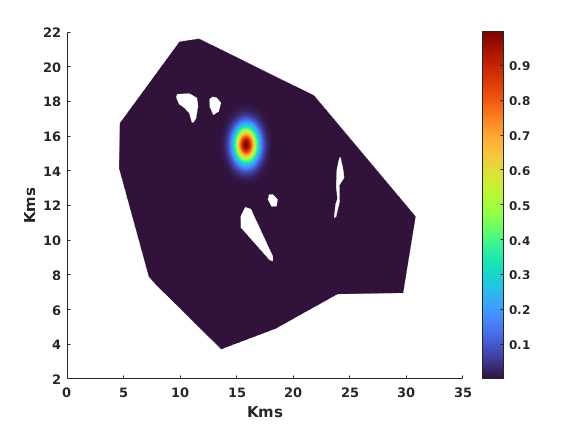}}
\subcaptionbox{Initial density $\rho^0$\label{fig3:initialdensity}}{\includegraphics[width=.8\linewidth,height=.2\textheight,keepaspectratio]{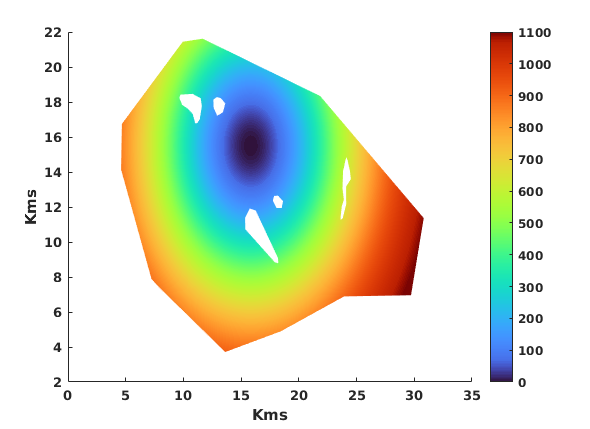}}
\subcaptionbox{Eikonal potential $\phi$\label{fig3:eikonal}}{\includegraphics[width=.8\linewidth,height=.2\textheight,keepaspectratio]{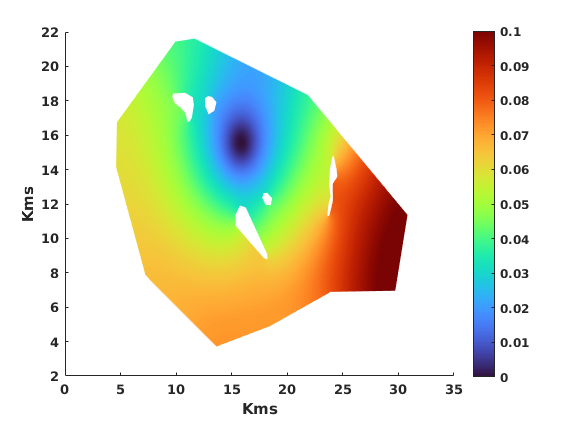}}
\subcaptionbox{Desired speed direction at initial time $\mathbf{v}(\rho^0)$\label{fig3:desiredspeedvector}}{\includegraphics[width=.8\linewidth,height=.2\textheight,keepaspectratio]{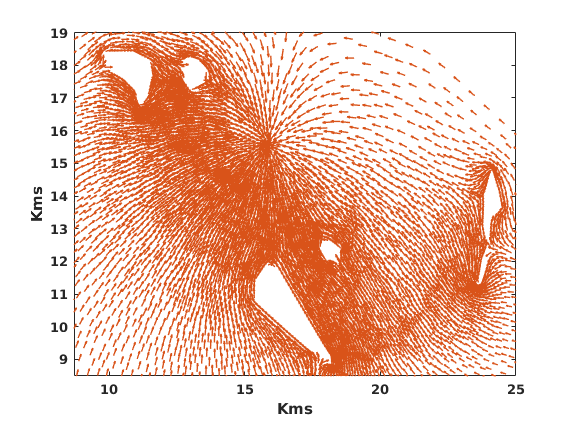}}
\subcaptionbox{Magnitude of initial desired speed $\Vert\mathbf{v}(\rho^0)\Vert$\label{fig3:desiredspeedvelocity}}{\includegraphics[width=.8\linewidth,height=.2\textheight,keepaspectratio]{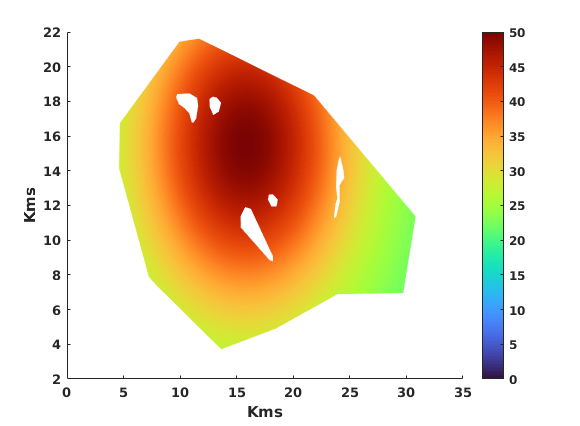}}
\subcaptionbox{Absorption rate of the solid phase $\epsilon \kappa$\label{fig3:kappa}}{\includegraphics[width=.8\linewidth,height=.2\textheight,keepaspectratio]{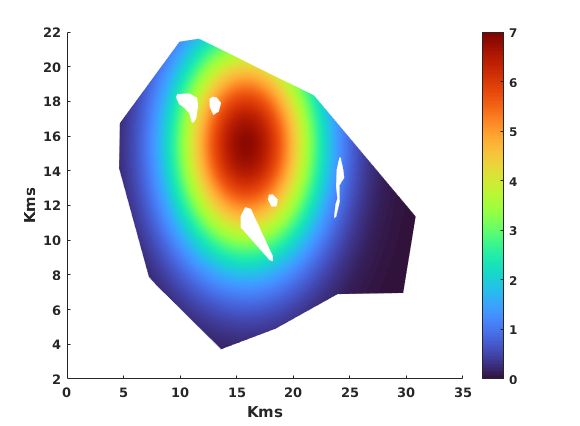}}
\caption{Inputs and initial state of the traffic density for our model: (\subref{fig3:G}) Normalized function $G$ used to enforce Eikonal potential $\phi$ to target the city center.  (\subref{fig3:initialdensity}) Initial state of the traffic density $\rho^0$, constructed under the idea of a concentric city.  (\subref{fig3:eikonal}) Eikonal potential $\phi$ derived from the definition of $G$.  (\subref{fig3:desiredspeedvector}) Detail of the vector field for the desired speed $\mathbf{v}(\rho^0)$,  pointing to city center and avoiding obstacles.  (\subref{fig3:desiredspeedvelocity}) Euclidean norm of $\mathbf{v}(\rho^0)$,  reaching -but not surpassing- maximum velocity $U_{max}=50$.   (\subref{fig3:kappa}) Absorption rate distribution in the whole city,  assuming -as in previous parameters- a concentric structure of the city.}\label{fig2}
\end{figure}

Meanwhile,  some functions defined on all nodes in the mesh $\mathcal{K}_h$ are needed to complement the model. The function $G$ introduced on the modified equation (\ref{ModPsi2:eq}) is displayed in Figure \ref{fig3:G},  clearly showing where the attraction point is. The initial state of the traffic density $\rho^0$ was constructed with the idea of the concentric city, with a density of $1000\ {\rm veh/Km}^2$ on the environs,  but less than $100$ in the downtown (Figure \ref{fig3:initialdensity}).  From this initial state, we show in Figures \ref{fig3:eikonal}, \ref{fig3:desiredspeedvector} and \ref{fig3:desiredspeedvelocity} the distribution of the travel cost $\phi$,  the desired speed direction $\mathbf{v}(\rho^0)$,  and its velocity values $\Vert\mathbf{v}(\rho^0)\Vert$,  respectively.  We can see here how all fields are consistent with the function $G$, since $\phi$ presents also a concentric distribution,  and $\mathbf{v}(\rho^0)$ targets the attraction point,  reaching maximum velocities ($U_{max}$) due its independence on porosity.  Finally, the capacity of drivers to find a parking space (i.e., the absorption rate $\kappa$) is also defined with a Gaussian distribution centered on the attraction point with a maximum absorption of $\epsilon \kappa =7$, and almost zero at the city limits (see Figure \ref{fig3:kappa}).

Although we have developed many numerical experiments to demonstrate the capabilities of our model, and the influence of its key parameters,  for the sake of conciseness we will present here only a few ones.  We will begin by showing the influence of the urban porosity on the behaviour of the model. 

\subsection{Two city scenarios modeled by the urban porosity: the dense city and the disperse city}

\begin{figure}
\centering
\subcaptionbox{The dense city\label{fig4:CiudadDensa}}{\includegraphics[width=.8\linewidth,height=.2\textheight,keepaspectratio]{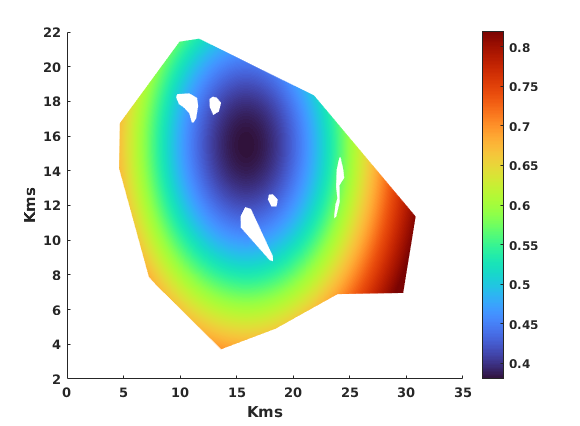}}
\subcaptionbox{The disperse city\label{fig4:CiudadDispersa}}{\includegraphics[width=.8\linewidth,height=.2\textheight,keepaspectratio]{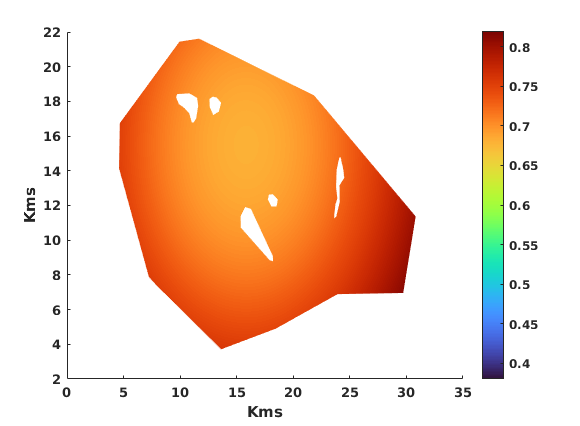}}
\caption{The two city scenarios are defined by their urban porosity values: (\subref{fig4:CiudadDensa}) The dense building city presents a porosity value at center $\epsilon_c = 0.38$.  (\subref{fig4:CiudadDispersa}) The disperse building city shows a porosity value at the city center $\epsilon_c=0.62$.  In both cases, the city environs have a porosity maximum value of $\epsilon_{max} = 0.82$.}\label{fig4}
\end{figure}

\begin{figure}
\centering
\subcaptionbox{Traffic density for dense city\label{fig5:densidadciudaddensa}}{\includegraphics[width=.8\linewidth,height=.2\textheight,keepaspectratio]{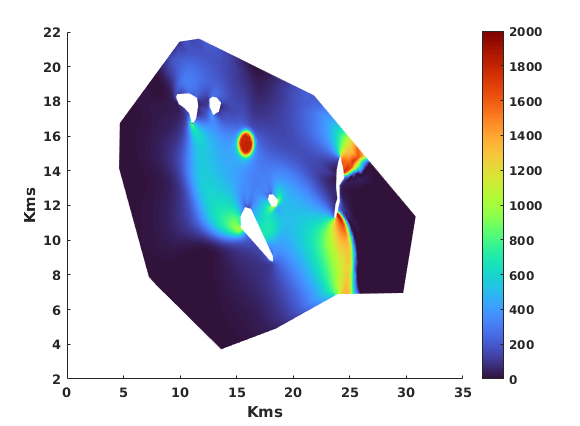}}
\subcaptionbox{Traffic density for disperse city\label{fig5:densidadciudaddispersa}}{\includegraphics[width=.8\linewidth,height=.2\textheight,keepaspectratio]{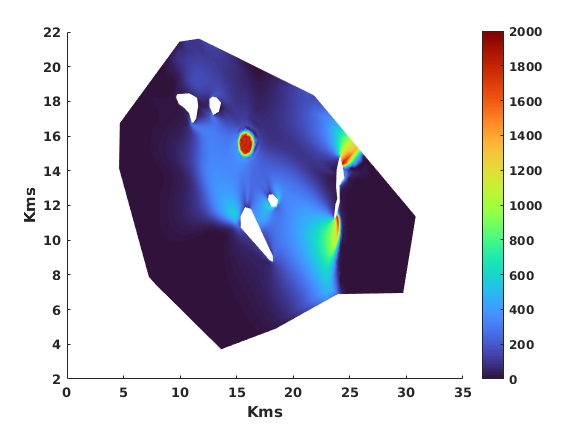}}
\subcaptionbox{Traffic speed for dense city\label{fig5:Unormaciudaddensa}}{\includegraphics[width=.8\linewidth,height=.2\textheight,keepaspectratio]{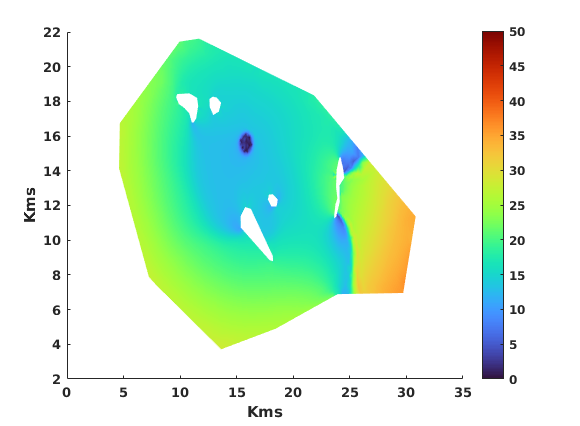}}
\subcaptionbox{Traffic speed for disperse city\label{fig5:Unormaciudaddispersa}}{\includegraphics[width=.8\linewidth,height=.2\textheight,keepaspectratio]{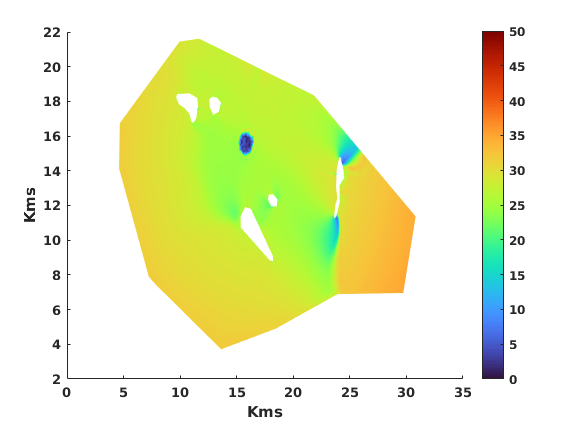}}
\caption{Traffic density and speed for the two city scenarios at the same simulation time $t=0.25$ h.  In  (\subref{fig5:densidadciudaddensa}) and (\subref{fig5:Unormaciudaddensa}) we can observe how, if buildings cover more area of the city, drivers must advance slow, reaching only low velocities, and taking more time to approach the city center after gathering around obstacles.  Unlike this,  as shown in (\subref{fig5:densidadciudaddispersa}) and (\subref{fig5:Unormaciudaddispersa}), an area with more streets facilitates the traffic flow, allowing higher velocities, and spending less time to reach the city center.  \label{fig5}}
\end{figure}

Following the interpretation of the porosity as a building density indicator, we define two very different city scenarios: a dense city with porosity values between $\epsilon_c = 0.38$ at downtown and $\epsilon_{max} = 0.82$ at city limits, and a disperse city with the same value at boundary $\epsilon_{max}=0.82$ but $\epsilon_c=0.62$ at city center (see Figure \ref{fig4}). Thus, the dense city has most of its area covered by buildings and less covered by streets,  being the opposite for the disperse city, with less area covered by buildings and more area covered by streets.  As it seems obvious, with more or less area available for the traffic flow, we expect an observable impact on traffic speed and traffic jams. The results obtained after $0.25$ hours of simulation can be seen in Figure \ref{fig5}: the dense city has more congested areas near the obstacles where the traffic velocity reaches only between $0$ and $15$ Km/h. In the case of the disperse city, traffic reaches a higher speed in the same areas (between $0$ and $25$ Km/h).  Therefore, in this case, cars can reach the attraction point faster (qualitative and quantitative differences between both scenarios can be seen in Figure \ref{fig5}). However, we can observe how, in both cases, the city center presents congestions due to a low interchange rate at downtown ($\kappa_{max} = 18$).  

\subsection{Influence of some key parameters of the model: absorption rate $\kappa$ and relaxation time $\tau$}

\begin{figure}
\centering
\subcaptionbox{Traffic density,  $\tau = 0.09$\label{figtau:density09}}{\includegraphics[width=.8\linewidth,height=.2\textheight,keepaspectratio]{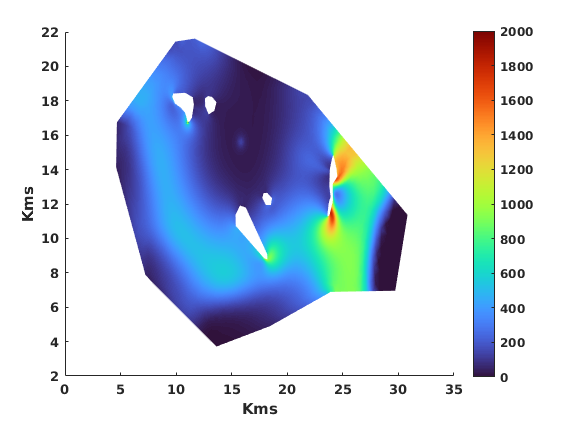}}
\subcaptionbox{Traffic speed,  $\tau = 0.09$\label{figtau:speed09}}{\includegraphics[width=.8\linewidth,height=.2\textheight,keepaspectratio]{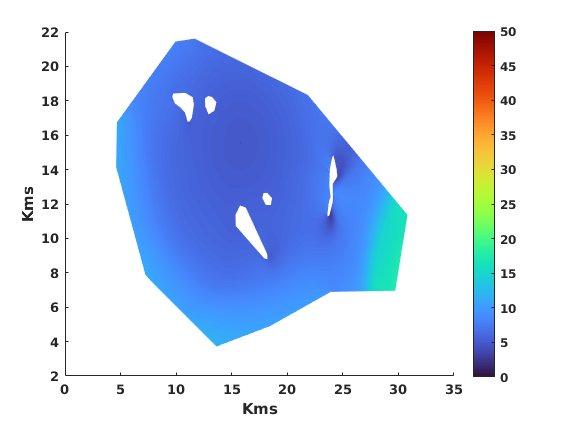}}
\caption{Influence of a larger relaxation time $\tau$ on the traffic density and speed: For $\tau = 0.09$, drivers require more time to approach the desired speed. This results in a low traffic speed (\subref{figtau:speed09}) with cars far to achieve the city center (\subref{figtau:density09}).  These pictures must be compared to Figure \ref{fig5}, corresponding to the smaller relaxation time $\tau=0.009$,  where the traffic speed approaches to the desired speed involving higher speed values (\subref{fig5:Unormaciudaddensa}), allowing cars to overrun obstacles, and congesting the city center (\subref{fig5:densidadciudaddensa}).  \label{figtau}}
\end{figure}

From all the parameters influencing the model we highlight two very important ones: the relaxation time $\tau$, and the absorption rate $\kappa$.  Relaxation time is defined as the time that drivers delay reaching the desired speed $\mathbf{v}(\rho)$, meanwhile, the absorption rate is related to the number of vehicles that pass from streets (fluid-phase) to parking spaces inside buildings or parking lines on streets (solid-phase). It is important remarking here that the desired speed depends on density but not on porosity.  Thus, the desired speed is an ideal speed where drivers could flow in a city without buildings reaching the higher possible velocities. Therefore,  we expect that, if $\tau \rightarrow 0$,  then $\mathbf{u}\rightarrow\mathbf{v}(\rho)$, neglecting the influence of the urban porosity.  In contrast to this, if $\tau\rightarrow 1$ we expect a full influence of porosity and, consequently, lower velocities. Figure \ref{figtau} shows, at the time $0.25$ h, the traffic density, and speed for $\tau=0.09$ in the dense city,  that must be compared to Figures \ref{fig5:densidadciudaddensa} and \ref{fig5:Unormaciudaddensa} corresponding to $\tau=0.009$.  In this latter case,  $\tau=0.009$ indicates a lower time of relaxation, and  consequently the traffic density presents a more congested downtown: is in the middle of overrunning the large obstacle in the right zone of the city,  having almost surpassing the other obstacles. This is due to velocities between $0$ and $30$ Km/h in the most congested zones.  Opposite to this, for a larger relaxation time $\tau=0.09$,  the cars are still arriving to obstacles despite the same time of simulation.  This is a consequence of lower velocity values in all the city, due to a higher influence of porosity, leaving a less congested center since a lower number of cars has reached it.

\begin{figure}
\centering
\subcaptionbox{Traffic Density \label{figkappa:density}}{\includegraphics[width=.8\linewidth,height=.2\textheight,keepaspectratio]{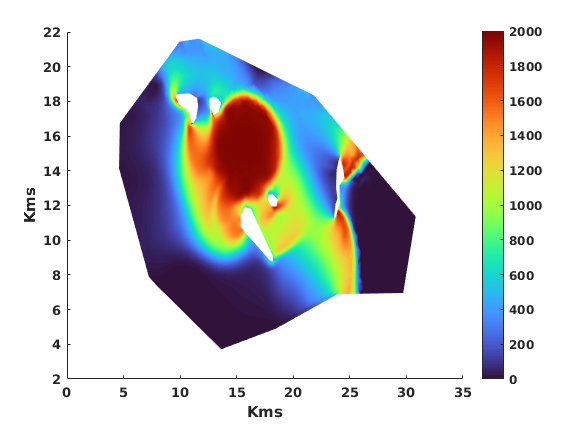}}
\subcaptionbox{Traffic Speed \label{figkappa:speed}}{\includegraphics[width=.8\linewidth,height=.2\textheight,keepaspectratio]{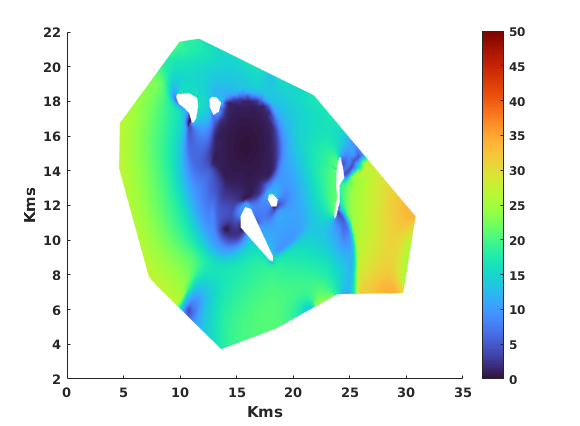}}
\caption{Influence of a lower absorption rate $\kappa$ on the traffic density and speed: Results obtained,  with $\kappa=1.8$,  for traffic density (\subref{figkappa:density}) and  traffic speed (\subref{figkappa:speed}).  A lower absorption rate implies more cars remaining in the streets,  blocking them,  increasing the congestion, and reducing the velocity to zero.\label{figkappa}}
\end{figure}

With respect to the role of the absorption rate in the model, we recall that it is a function with a higher value at the city center $\kappa_{max}$ tending to zero as we reach the city limits (defined again as a Gaussian distribution). Thus, since we already used the value $\kappa_{max}=18$ in the previous experiment,  now we present here the traffic density and speed corresponding to a lower value of the absorption value ($\kappa_{max} = 1.8$). So, in Figures \ref{figkappa:density} and \ref{figkappa:speed} we can see the traffic density and speed, respectively, for this new value.  As expected, a lower absorption rate can be interpreted as the fact that drivers cannot find a parking space, remaining on the streets and, as more cars reach the center, they are forced to occupy a bigger area (this is the wide red area observed in Figure \ref{figkappa:density}) with almost zero velocities values (the large blue area in Figure \ref{figkappa:speed}).

\subsection{The traffic demand $(1-\epsilon)q$}

The traffic demand represents the traffic source in the model.  In our interpretation of the traffic flow at an urban porous media, the traffic demand simulates the flow that brings cars from buildings (solid-phase) to streets (fluid-phase). Also, to be consistent with a concentric type of city,  its maximum values must be located in residential zones near the city limits, with minimum values at the city downtown.  Moreover, traffic demand must consider the local travel cost,  which in our case involves a zero cost if cars are far from the center but increases conforming cars reach it \cite{Jiang2016}.  Finally,  traffic demand must also incorporate a time component to simulate rush-valley hours along the day.  Therefore, since traffic demand must depend on the solid-phase $(1-\epsilon)$, on the local travel cost $(1-\phi_{max}/\phi)$,  on a suitable function $q_{max}(\mathbf{x})$ defined for all $\mathbf{x}\in {\Omega}$, and on a time function $g(t)$,  we are led to formulate the traffic demand function as:
$$
q = q_{max}\big(1-\frac{\phi_{max}}{\phi}\big)\,g(t).
$$

\begin{figure}
\centering
\subcaptionbox{Travel local cost $(1-\phi_{max}/\phi)$\label{fig8q:travelcost}}{\includegraphics[width=.8\linewidth,height=.2\textheight,keepaspectratio]{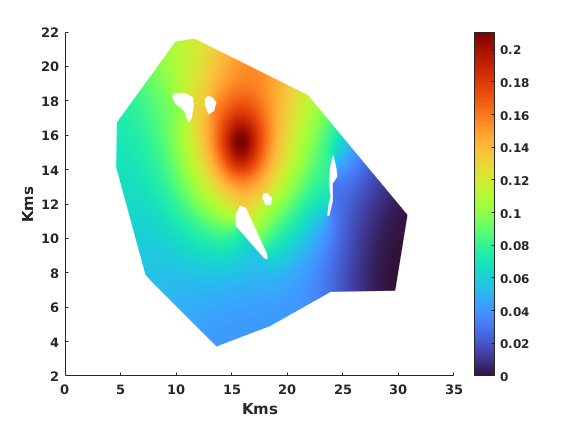}}
\subcaptionbox{Maximum traffic demand from the solid-phase $(1-\epsilon)q_{max}$\label{fig8q:trafficdemand}}{\includegraphics[width=.8\linewidth,height=.2\textheight,keepaspectratio]{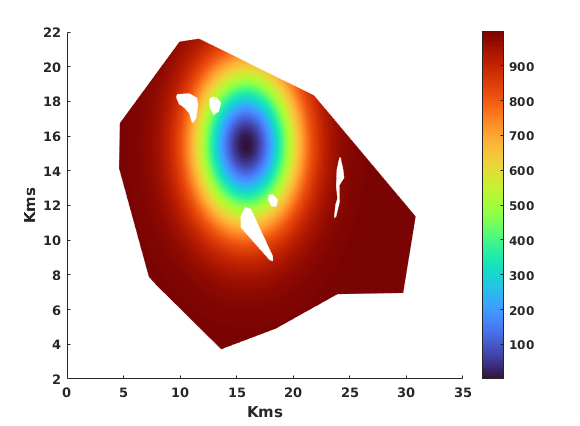}}
\subcaptionbox{The $g(t)$ function during the $4$ hours of simulation time \label{fig8q:gt}}{\includegraphics[width=.8\linewidth,height=.2\textheight,keepaspectratio]{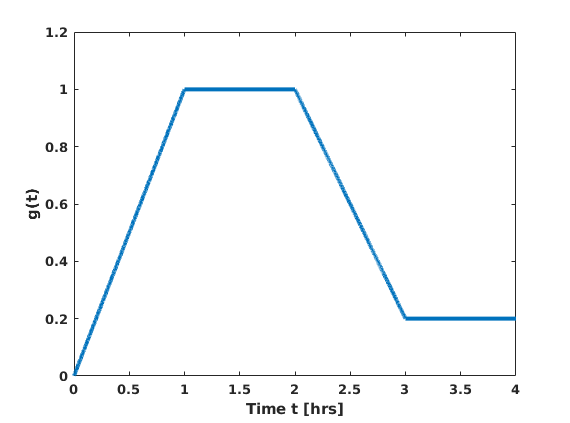}}
\caption{The local travel cost $(1-\phi_{max}/\phi)$ at the first instant of simulation $t=0$,  the spatial distribution of traffic demand $q_{max}$ from buildings (solid-phase), and the time distribution $g(t)$ showing a rush-valley hours behaviour.}\label{fig8q}
\end{figure}

\begin{figure}\captionsetup[figure]{font=footnotesize}
\centering
\subcaptionbox{$t=0.5$ h \label{fig9q:mediahora}}{\includegraphics[width=.8\linewidth,height=.2\textheight,keepaspectratio]{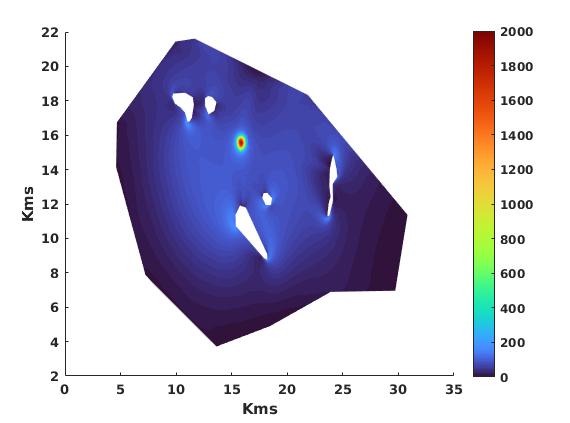}}
\subcaptionbox{$t=1.5$ h \label{fig9q:1horasymedia}}{\includegraphics[width=.8\linewidth,height=.2\textheight,keepaspectratio]{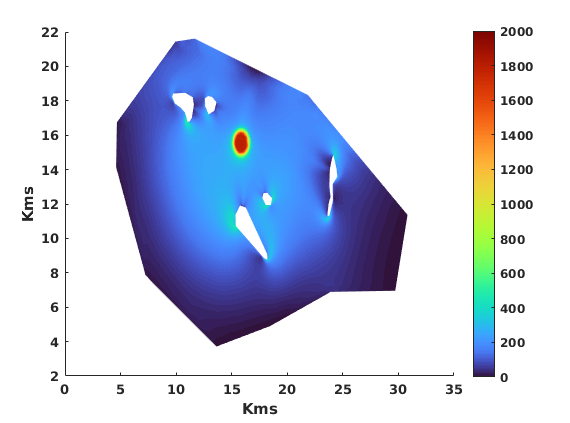}}
\subcaptionbox{$t=2.5$ h \label{fig9q:2horasymedia}}{\includegraphics[width=.8\linewidth,height=.2\textheight,keepaspectratio]{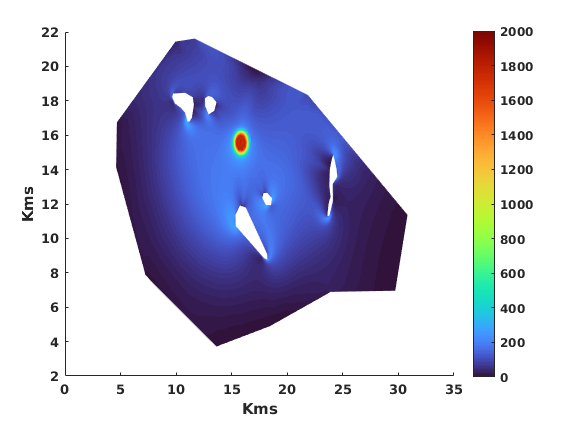}}
\subcaptionbox{$t=3.5$ h \label{fig9q:3horasymedia}}{\includegraphics[width=.8\linewidth,height=.2\textheight,keepaspectratio]{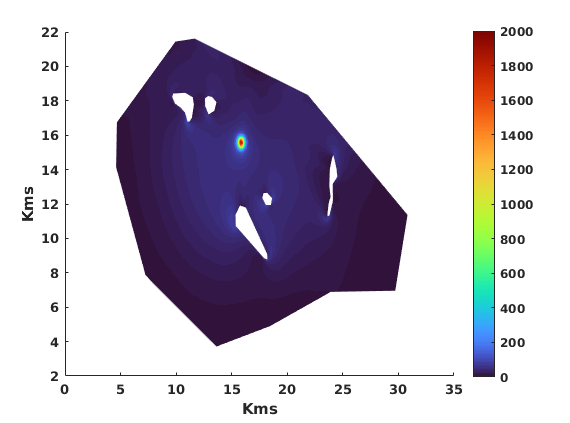}}
\caption{Traffic density at different stages in a typical traffic period passing from the first hour increase (\subref{fig9q:mediahora}) to a rush hour (\subref{fig9q:1horasymedia}), where density reaches its peak with a congested center.  After the third hour (\subref{fig9q:2horasymedia}) the traffic demand decreases,  and so density,  since most of the drivers have already reached their destination (\subref{fig9q:3horasymedia}).} \label{fig9q}
\end{figure}

In Figures \ref{fig8q:travelcost} and \ref{fig8q:trafficdemand} we display the travel cost $(1-\phi_{max}/\phi)$ and the maximum traffic demand $(1-\epsilon)q_{max}$ at the beginning of the simulation. The travel cost is zero for cars located in the remote part of the city with respect to the city center, the cost increases as cars address the attraction point.  Meanwhile, the traffic demand shows a concentric distribution,  achieving its maximum values on a bell surrounding the center. The traffic demand is controlled by function $g(t)$, which takes values between $[0,1]$, being $g=1$ a rush hour and $g=0.2$ a valley hour (see Figure \ref{fig8q:gt}).  Once the traffic demand was defined,  we ran an experiment with a time simulation of 4 hours,  passing by a demand increase in the first hour (this is, $g(t)$ increases from $0$ to $1$), a maximum demand during a rush hour ($g(t) = 1$),  and finally a traffic decline where $g(t)$ passes from $1$ to $0.2$,  remaining in this latter value until simulation time ends (see Figure \ref{fig8q:gt}).  Recalling that drivers travel towards the attraction point at the city center, we can see that, as the traffic demand increases in time, some cars may leave the streets in their path, but most of them reach the center congesting it (Figures \ref{fig9q:mediahora} and \ref{fig9q:1horasymedia}).  Then, as the demand decreases, a lower number of cars take the path to the center and, consequently, the cars can find easily a parking space (Figure \ref{fig9q:2horasymedia}),  allowing the center to recover its empty streets (Figure \ref{fig9q:3horasymedia}).

\section{Conclusions}\label{sec_conclusions}




\begin{thebibliography}{10}
\expandafter\ifx\csname url\endcsname\relax
  \def\url#1{\texttt{#1}}\fi
\expandafter\ifx\csname urlprefix\endcsname\relax\def\urlprefix{URL }\fi
\expandafter\ifx\csname href\endcsname\relax
  \def\href#1#2{#2} \def\path#1{#1}\fi

\bibitem{hu2012}
Z.~Hu, B.~Yu, Z.~Chen, T.~Li, M.~Liu, Numerical investigation on the urban heat
  island in an entire city with an urban porous media model, Atmospheric
  Environment 47 (2012) 509--518.
\newblock \href {https://doi.org/10.1016/j.atmosenv.2011.09.064}
  {\path{https://doi.org/10.1016/j.atmosenv.2011.09.064}}.

\bibitem{wang2021}
H.~Wang, C.~Peng, W.~Li, C.~Ding, T.~Ming, N.~Zhou,
  \href{https://www.sciencedirect.com/science/article/pii/S221209552100095X}{Porous
  media: A faster numerical simulation method applicable to real urban
  communities}, Urban Climate 38 (2021) 100865.
\newblock \href {https://doi.org/https://doi.org/10.1016/j.uclim.2021.100865}
  {\path{https://doi.org/https://doi.org/10.1016/j.uclim.2021.100865}}.

\bibitem{ming2021}
T.~Ming, S.~Lian, Y.~Wu, T.~Shi, C.~Peng, Y.~Fang, R.~de~Richter, N.~H. Wong,
  \href{https://www.mdpi.com/1996-1073/14/15/4681}{Numerical investigation on
  the urban heat island effect by using a porous media model}, Energies 14~(15) 
  (2021) 4681.
\newblock \href {https://doi.org/10.3390/en14154681}
  {\path{https://doi.org/10.3390/en14154681}}.

\bibitem{garcia-chan2023}
N.~García-Chan, J.~A. Licea-Salazar, L.~G. Gutierrez-Ibarra,
  \href{https://www.mdpi.com/2227-7390/11/5/1140}{Urban heat island dynamics in
  an urban--rural domain with variable porosity: Numerical methodology and
  simulation}, Mathematics 11~(5) (2023) 1140.
\newblock \href {https://doi.org/10.3390/math11051140}
  {\path{https://doi.org/10.3390/math11051140}}.

\bibitem{treiber2013}
A.~K. Martin~Treiber, Traffic Flow Dynamics, Springer, Berlin, 2013.

\bibitem{Jiang2014}
Y.~Q. Jiang, S.~G. Zhou, Macroscopic simulation of traffic flow on continuum
  urban networks, Applied Mechanics and Materials (2014) 641--642.  \href
  {https://doi.org/10.4028/www.scientific.net/amm.641-642.887}
  {\path{https://doi:10.4028/www.scientific.net/amm.641-642.887}}.

\bibitem{Jiang2016}
Y.~Q. Jiang, P.~J. Ma, S.~G. Zhou, Macroscopic modeling approach to estimate
  traffic-related emissions in urban areas, Transportation Research Part D:
  Transport and Environment 60 (2018) 41--55.
\newblock \href {https://doi.org/10.1016/j.trd.2015.10.022}
  {\path{https://doi.org/10.1016/j.trd.2015.10.022}}

\bibitem{huber2015}
G.~Huber, S.~Tanguy, J.-C. Béra, B.~Gilles,
  \href{https://www.sciencedirect.com/science/article/pii/S0021999115006129}{A
  time splitting projection scheme for compressible two-phase flows:
  application to the interaction of bubbles with ultrasound waves}, Journal of
  Computational Physics 302 (2015) 439--468.
\newblock \href {https://doi.org/10.1016/j.jcp.2015.09.019}
  {\path{https://doi.org/10.1016/j.jcp.2015.09.019}}.

\bibitem{tezduyar2006}
T.~E. Tezduyar, M.~Senga,
  \href{https://www.sciencedirect.com/science/article/pii/S0045782505002999}{Stabilization
  and shock-capturing parameters in supg formulation of compressible flows},
  Computer Methods in Applied Mechanics and Engineering 195~(13) (2006)
  1621--1632.
\newblock \href {https://doi.org/10.1016/j.cma.2005.05.032}
  {\path{https://doi.org/10.1016/j.cma.2005.05.032}}.

\bibitem{caltagirone2011}
J.~P. Caltagirone, S.~Vincent, C.~Caruyer,
  \href{https://www.sciencedirect.com/science/article/pii/S004579301100199X}{A
  multiphase compressible model for the simulation of multiphase flows},
  Computers \& Fluids 50~(1) (2011) 24--34.
\newblock \href
  {https://doi.org/10.1016/j.compfluid.2011.06.011}
  {\path{https://doi.org/10.1016/j.compfluid.2011.06.011}}.

\bibitem{payne1971}
H.~Payne, \href{https://books.google.com.mx/books?id=1X4ZcgAACAAJ}{Models of
  Freeway Traffic and Control}, Simulation Councils, Inc., 1971.

\bibitem{kerner_konhauser1993}
B.~S. Kerner, P.~Konh\"auser,
  \href{https://link.aps.org/doi/10.1103/PhysRevE.48.R2335}{Cluster effect in
  initially homogeneous traffic flow}, Phys. Rev. E 48 (1993) R2335--R2338.
\newblock \href {https://doi.org/10.1103/PhysRevE.48.R2335}
  {\path{https://doi.org/10.1103/PhysRevE.48.R2335}}.

\bibitem{herman-tenny-prigrione1971}
R.~Herman, T.~Lam, I.~Prigogine,
  \href{http://www.jstor.org/stable/25767676}{Kinetic theory of vehicular
  traffic: Comparison with data}, Transportation Science 6~(4) (1972) 440--452.
\newblock \href {https://doi.org/10.1287/trsc.6.4.440}
  {\path{https://doi.org/10.1287/trsc.6.4.440}}

\bibitem{daganzo1995}
C.~F. Daganzo,
  \href{https://www.sciencedirect.com/science/article/pii/019126159500007Z}{Requiem
  for second-order fluid approximations of traffic flow}, Transportation
  Research Part B: Methodological 29~(4) (1995) 277--286.
\newblock \href {https://doi.org/10.1016/0191-2615(95)00007-Z}
  {\path{https://doi.org/10.1016/0191-2615(95)00007-Z}}.

\bibitem{Aw2000}
A.~Aw, M.~Rascle, \href{https://doi.org/10.1137/S0036139997332099}{Resurrection
  of ``second order" models of traffic flow}, SIAM Journal on Applied
  Mathematics 60~(3) (2000) 916--938.
\href  {https://doi.org/10.1137/S0036139997332099}
  {\path{https://doi.org/10.1137/S0036139997332099}}.

\bibitem{das2018_book}
M.~K. Das, P.~P. Mukherjee, K.~Muralidhar, Modeling Transport Phenomena in
  Porous Media with Applications, Springer International Publishing, Cham,  2018.

\bibitem{lsm2017}
L.~Leclercq,  A.~Senecata,  G.~Mariotte, 
\href{https://doi.org/10.1016/j.trb.2017.04.004}{Dynamic macroscopic 
simulation of on-street parking search: A trip-based approach}, Transportation Research 
Part B 101 (2017) 268--282.
\href  {https://doi.org/10.1016/j.trb.2017.04.004}
  {\path{https://doi.org/10.1016/j.trb.2017.04.004}}.

\bibitem{sportisse2010}
B.~Sportisse, Fundamentals in Air Pollution,  Springer Netherlands,
  2010.

\bibitem{Xia2008}
Y.~Xia, S.~C. Wong, M.~Zhang, C.-W. Shu, W.~H.~K. Lam, An efficient
  discontinuous galerkin method on triangular meshes for a pedestrian flow
  model, International Journal for Numerical Methods in Engineering 76~(3)
  (2008) 337--350.
\newblock \href {https://doi.org/https://doi.org/10.1002/nme.2329}
  {\path{https://doi.org/10.1002/nme.2329}}.

\bibitem{Jiang2015}
Y.~Q. Jiang, S.~G. Zhou, F.~B. Tian, A higher-order macroscopic model for
  bi-direction pedestrian flow, Physica A: Statistical Mechanics and its
  Applications 425 (2015) 69--78.
\newblock \href {/https://doi.org/10.1016/j.physa.2014.11.048}
  {\path{https://doi.org/10.1016/j.physa.2014.11.048}}.

\bibitem{axthelm2016}
R.~Axthelm, Finite element simulation of a macroscopic model for pedestrian
  flow, in: V.~L. Knoop, W.~Daamen (Eds.), Traffic and Granular Flow'15,
  Springer International Publishing, Cham, 2016, pp. 233--240.

\bibitem{larson}
M.~G. Larson, F.~Bengzon, The Finite Element Method: Theory, Implementation,
  and Applications, Springer,  Berlin, 2013.

\bibitem{gottlieb_ketcheson2009}
S.~Gottlieb, D.~Ketcheson, C.~Shu, High order strong stability preserving time
  discretizations, Journal of Scientific Computing 38~(3) (2009) 251--289.
\newblock \href {https://doi.org/10.1007/s10915-008-9239-z}
  {\path{https://doi.org/10.1007/s10915-008-9239-z}}.

\bibitem{Knaber2003}
P.~Knaber, L.~Angermann, Numerical methods for Elliptic and Parabolic Partial Differential Equations,
Springer-Verlag,  New York, 2003.

\bibitem{Ern2021}
A.~Ern, J. L.~Guermond, Finite Elements III, First-Order and Time-Dependent PDEs, Springer Nature, Switzerland, 2021.

\bibitem{gmsh2009}
C.~Geuzaine, J.~F. Remacle,
  \href{https://onlinelibrary.wiley.com/doi/abs/10.1002/nme.2579}{Gmsh: A 3-D
  finite element mesh generator with built-in pre- and post-processing
  facilities}, International Journal for Numerical Methods in Engineering
  79~(11) (2009) 1309--1331.
\newblock   \href {https://doi.org/10.1002/nme.2579}
  {\path{https://doi.org/10.1002/nme.2579}}.

\end{thebibliography}


In this work, we introduce a novel macroscopic traffic flow model on urban domains with the main novelties being the city interpretation as a porous media and the non-conservation of cars. Both elements have not been considered in previous similar models appearing in the mathematical literature on the topic.  Indeed, the urban porosity was used to describe the urban landscape (considering the streets as the fluid-phase of a porous media, and the buildings as its solid-phase),  and the non-conservation allows us to model the typical traffic dynamics of cars leaving its garages to flow towards the city and make off-street parking when arriving to their destiny.

Concerning the numerical techniques used here, we highlight that the classical Galerkin $P_1$ finite element method combined with an explicit marching time scheme was enough to achieve stable numerical solutions.  This was not possible in previous works,  where hyperbolic systems of PDE were used.  The stability is a consequence of the density diffusion term included in the continuity equation, and of the Brinkman correction in the case of the Darcy Law,  which adds viscosity to the momentum equation. The stability of our algorithm allowed us to make several experiments testing the response of the model to the flow around obstacles, absorption of vehicles in the city center -off-street parking-, changes on the so-called relaxation time, changes in traffic demand, and a variable urban-porosity to characterize the dense and disperse cities (porosity gradients).

Our numerical results make evident the influence of the city landscape (porosity) on the traffic flow,  being faster in cities with large streets.  Contrary to this,  the flow is slower when more area is dedicated to building blocks.  Also, the absorption rate (parking capacity) is important to take out the great number of cars avoiding congestion.  On the opposite side,  a poor capacity makes that the congestion area covers larger city zones.  The model response to limited relaxation time was to neglect the urban porosity,  because the Eikonal model did not consider it.  Meanwhile, the transportation cost warrants a car's density below its maximum value.  Even though these results are expected in a concentric city being its interpretation straightforward, our methodology holds equal for more complicated cities with irregular limit layouts as well as a chaotic porosity distribution.

Promising results obtained here,  allow us to move on to more complex situations, such as the existence of different parking spaces (i.e., multiple attraction points), or the incorporation of other phenomena not considered in this work.





\end{document}